\documentclass[final,3p,times]{elsarticle}

\usepackage{amssymb}

\usepackage{amsmath}
\usepackage{multirow}
\usepackage{color}

\journal{Journal of Computational Physics}

\begin{document}

\begin{frontmatter}



\title{High-Order Schemes of Exponential Time Differencing for Stiff Systems\\ with Nondiagonal Linear Part}


\author[icmm]{Evelina V. Permyakova}
\ead{evelina.v.permyakova@gmail.com}

\author[icmm,psu]{Denis S. Goldobin\corref{cor1}}
\ead{Denis.Goldobin@gmail.com}

\cortext[cor1]{Corresponding author}
\address[icmm]{Institute of Continuous Media Mechanics UB RAS,
             1 Akad.\ Koroleva street,
             Perm 614013,
             Russia}

\address[psu]{Institute of Physics and Mathematics, Perm State University,
             15 Bukireva street,
             Perm 614990,
             Russia}

\begin{abstract}
Exponential time differencing  methods is a power tool for high-performance numerical simulation of computationally challenging problems in condensed matter physics, fluid dynamics, chemical and biological physics, where mathematical models often possess fast oscillating or decaying modes---in other words, are stiff systems. Practical implementation of these methods for the systems with nondiagonal linear part of equations is exacerbated by infeasibility of an analytical calculation of the exponential of a nondiagonal linear operator; 
in this case, the coefficients of the exponential time differencing scheme cannot be calculated analytically.  We suggest an approach, where these coefficients are numerically calculated with auxiliary problems. We rewrite the high-order Runge--Kutta type schemes in terms of the solutions to these auxiliary problems and practically examine the accuracy and computational performance of these methods for a heterogeneous Cahn--Hilliard equation, a sixth-order spatial derivative equation governing pattern formation in the presence of an additional conservation law, and a Fokker--Planck equation governing macroscopic dynamics of a network of neurons.
\end{abstract}


\begin{highlights}
\item An approach to implementation of exponential time differencing methods is presented
\item The approach allows for simulation of equation systems with nondiagonal linear part
\item Performance gain for high-precision calculations is up to several orders of magnitude
\item Efficient for PDEs, Fokker-Planck equations, and large ensembles of oscillators
\end{highlights}


\end{frontmatter}


\section{Introduction}
In condensed matter, chemical and biological physics, many systems are governed by the mathematical models which can be called `stiff' systems---the systems where some modes are fast oscillating or decaying. The dynamics of wave function in atom lattices with frozen parametric disorder in potential is fundamentally contributed by the fast oscillating modes and this dynamics gives rise to the phenomenon of Anderson localization~\cite{Anderson-1958,Frohlich-Spenser-1984,Lifshitz-Gredeskul-Pastur-1988,Blumel-Fishman-Smilansky-1985,Schwartz-etal-2007,Pikovsky-Shepelyansky-2008}. The spinodal decomposition~\cite{Cahn-Hilliard-1958,Golovin-etal-2001,Podolny-etal-2005,Watson-etal-2003,Speck-etal-2015,Kuramoto-Tsuzuki-1976} of two-component mixtures and heat-mass-transfer in active media~\cite{Knobloch-1990,Shtilman-Sivashinsky-1991,
Matthews-Cox-2000,Matthews-Cox-2000b,Mosheva-Siraev-Bratsun-2023,Goldobin-Shklyaeva-2008,Shklyaev-etal-2012,Samoilova-Nepomnyashchy-2019,Samoilova-Nepomnyashchy-2020,Kozlov-2023} is governed by partial differential equations with high-order spatial derivatives; these derivatives drive a fast decay of short wavelength perturbations. The systems near a SNIPER (saddle-node infinite period) bifurcation of chemical oscillations~\cite{Erban-etal-2009} are subject to thermal noise and can be subject to other noise sources; at the bifurcation point, these systems can be described by the same mathematical models as the quadratic integrate-and-fire neurons with noise and synaptic coupling~\cite{Ratas-Ryragas-2019,Goldobin-diVolo-Torcini-2021,diVolo-etal-2022,Zheng-Kotani-Jimbo-2021,Goldobin-Permyakova-Klimenko-2024,Pietras-Cestnik-Pikovsky-2023,Goldobin-etal-2023}. The numerical simulation of the Fokker--Planck equation of these chemical oscillation and neuronal systems is again the case of a stiff system.

Generally, fast oscillating or decaying modes impose severe limitation on the time stepsize, which has to be very small. For oscillating modes it is needed to maintain a demanded accuracy of simulation~\cite{Pikovsky-Shepelyansky-2008}. Fast decaying modes are dying-out and do not influence the dynamics of the chemical/physical system; therefore, there is no need for a high precision of their simulation. However, for a finite stepsize, fast decay can produce numerical instability if the stepsize in not small enough. As a result, the numerical simulation of stiff system requires the employment of {\it ad hoc} approaches and implicit schemes, which require a sophisticated individual mathematical preliminary work, is costly or infeasible in higher dimensions, and are helpless for fast oscillatory modes in conservative systems. Otherwise, the simulation of these problems should be performed with a tiny time stepsize and becomes extremely CPU (central processing unit) time consuming. An alternative approach is a relatively new class of methods of exponential time differencing (ETD)~\cite{Cox-Matthews-2002,Minchev-Wright-2005,Hochbruck-Ostermann-2010,Beylkin-Keiser-Vozovoi-1998,Holland-1994,Petropoulos-1997,Schuster-Christ-Fichtner-2000,Hochbruck-Lubich-Selhofer-1998,Tokman-2006,Kassam-Trefethen-2005}.

In ETD methods~\cite{Cox-Matthews-2002,Minchev-Wright-2005,Hochbruck-Ostermann-2010} the linear part of equations, which gives the origin to the fast modes, is solved exactly. For the nonlinear part, one constructs an exact solution for a quasistatic approximation of these terms during one time step in the case of a first-order scheme. For an $n$th-order scheme, one effectively constructs the exact solution for the case where the nonlinear part in one time step is an $(n-1)$-order polynomial of time~\cite{Cox-Matthews-2002}. This approach provides high resolution for temporal dynamics of fast modes and prevents numerical instabilities even for a quite large values of time stepsize.

The application of ETD methods is straightforward when the linear part is diagonal, which is abundant for spectral methods for the problems with periodic boundary conditions and Fourier basis. However, for spectral methods with other basis functions (e.g., Chebyshev polynomials) the linear part can become nondiagonal. Moreover, for problems with frozen heterogeneity of parameters~\cite{diVolo-etal-2022,Hammele-Schuler-Zimmermann-2006,Goldobin-Shklyaeva-2009,Goldobin-Shklyaeva-2013,Goldobin-2010,Goldobin-2019}, the linear part becomes essentially nondiagonal. In this case the analytical calculation of the coefficients of ETD schemes (considered in Sec.~\ref{sec:ETD}) can become not just difficult but completely impossible.

In this paper we develop an approach to calculation of the coefficients of high-order schemes of exponential time differencing for systems with a nondiagonal linear part.
In this approach the coefficients are calculated by means of numerical integration of auxiliary problems. This integration is carried-out with a plain standard method, like the predictor--corrector one, and a very small time stepsize $\tau_1$, ensuring numerical stability, but over a short time interval---the time step $\tau$ of the ETD scheme. The basic concept was suggested in~\cite{Permyakova-Goldobin-2020} and here we develop and apply the approach to the Runge--Kutta type schemes~\cite{Cox-Matthews-2002} of 3rd and 4th order. For the sake of illustration, we consider the numerical simulation of the discrete-space versions of partial differential equations: Cahn--Hilliard equation~\cite{Cahn-Hilliard-1958,Golovin-etal-2001} and Matthews--Cox equation~\cite{Matthews-Cox-2000,Matthews-Cox-2000b}, and the Fourier projection of a Fokker--Planck equation~\cite{Ratas-Ryragas-2019,Zheng-Kotani-Jimbo-2021,Goldobin-2021}. The performance gain varies from two to many orders of magnitude.
The technical implementation of the suggested approach and the bulk of the text below do not employ advanced topics from Linear algebra and Tensor calculus.

The paper is organized as follows. In Sec.~\ref{sec:ETD}, we provide three exponential differencing schemes of the Runge--Kutta type and rewrite them in terms of the solutions of the auxiliary problems. The procedure for constructing the solutions of the auxiliary problems is provided. In Sec.~\ref{sec:Examples}, we considers examples of the application of the ETD schemes to three generic mathematical models: Cahn--Hilliard equation (Sec.~\ref{sec:CHE}), sixth-order spatial derivative Matthews--Cox equation (Sec.~\ref{sec:MCE}), and Fokker--Planck equation for a network of theta-neurons (Sec.~\ref{sec:QIFs}). The accuracy and CPU time performance of the ETD schemes is analysed.
In Sec.~\ref{sec:Concl} we derive conclusions. At \ref{sec:app}, we examine the ``strong'' accuracy of the suggested approach to numerical solution of the auxiliary problems.

\section{\label{sec:ETD}Exponential time differencing}

We deal with the problem of numerical integration of the equations system of the following form:
\begin{equation}
\dot{\mathbf{u}}=\mathbf{L}\cdot\mathbf{u}+\mathbf{f}(\mathbf{u},t)\,,
\label{eq:ETD01}
\end{equation}
where $\mathbf{u}(t)$ is an $N$-component vector, $\mathbf{L}$ is an $[N\times N]$-matrix with time-independent elements, $\mathbf{f}(\mathbf{u},t)$ is the nonlinear part of equations. The decomposition of the equations into the linear part with time-independent coefficients, $\mathbf{L}\cdot\mathbf{u}$, and the nonlinear part is not unique and guided by the merit of convenience: $\mathbf{f}(\mathbf{u},t)$ can even contain a time-independent linear part. However, it is desirable that all the dominating sources of numerical instability should be collected into $\mathbf{L}\cdot\mathbf{u}$ as fully as possible.

In this paper we consider implementation of three exponential differencing schemes of the Runge--Kutta type\cite{Cox-Matthews-2002}; for equation system~(\ref{eq:ETD01}) these schemes take the following form:
\\
\textbf{(ETD2RK)}~two-step scheme:
\begin{subequations}
\label{eq:ETD02ab}
\begin{align}
\mathbf{a}&=e^{\mathbf{L}\tau}\cdot\mathbf{u}(t) +\mathbf{L}^{-1}\cdot(e^{\mathbf{L}\tau}-\mathbf{I})\cdot\mathbf{f}\big(\mathbf{u}(t),t\big)\,,
\label{eq:ETD02}
\\
\mathbf{u}(t+\tau)&=\mathbf{a} +\mathbf{L}^{-2}\cdot(e^{\mathbf{L}\tau}-\mathbf{I}-\mathbf{L}\tau)\cdot\frac{\mathbf{f}\big(\mathbf{a},t+\tau\big) -\mathbf{f}\big(\mathbf{u}(t),t\big)}{\tau}\,,
\label{eq:ETD03}
\end{align}
\end{subequations}
where $\tau$ is the time stepsize of the numerical scheme, $\mathbf{a}$ is a preliminary approximation of $\mathbf{u}$ for the next time step $t+\tau$, $\mathrm{I}$ is the unitary matrix, the exponential of a matrix is defined by the series $\exp{\mathbf{A}}=\mathbf{I}+\mathbf{A}+\mathbf{A}^2/2!+\mathbf{A}^3/3!+\dots$;
\\
\textbf{(ETD3RK)}~three-step scheme:
\begin{subequations}
\label{eq:ETD03ac}
\begin{align}
\mathbf{a}&=e^{\mathbf{L}\tau/2}\cdot\mathbf{u}(t) +\mathbf{L}^{-1}\cdot(e^{\mathbf{L}\tau/2}-\mathbf{I})\cdot\mathbf{f}\big(\mathbf{u}(t),t\big)\,,
\label{eq:ETD04}
\\
\mathbf{b}&=e^{\mathbf{L}\tau}\cdot\mathbf{u}(t) +\mathbf{L}^{-1}\cdot(e^{\mathbf{L}\tau}-\mathbf{I})\cdot \left[2\mathbf{f}\big(\mathbf{a},t+\tau/2\big)-\mathbf{f}\big(\mathbf{u}(t),t\big)\right]\,,
\label{eq:ETD05}
\\
\mathbf{u}(t+\tau)&=e^{\mathbf{L}\tau}\cdot\mathbf{u}(t) +\frac{\mathbf{L}^{-3}}{\tau^2}\cdot\Big(\big[-4-\mathbf{L}\tau
+e^{\mathbf{L}\tau}(4-3\mathbf{L}\tau+\mathbf{L}^2\tau^2)\big]\cdot\mathbf{f}\big(\mathbf{u}(t),t\big)
\nonumber\\
&\qquad{}
+4\left[2+\mathbf{L}\tau +e^{\mathbf{L}\tau}(-2+\mathbf{L}\tau)\right]\cdot\mathbf{f}\big(\mathbf{a},t+\tau/2\big)
+\left[-4-3\mathbf{L}\tau-\mathbf{L}^2\tau^2 +e^{\mathbf{L}\tau}(4-\mathbf{L}\tau)\right]\cdot\mathbf{f}\big(\mathbf{b},t+\tau\big)
\Big)\,,
\label{eq:ETD06}
\end{align}
\end{subequations}
where $\mathbf{a}$ and $\mathbf{b}$ are the preliminary approximations of $\mathbf{u}$ for $t+\tau/2$ and $t+\tau$, respectively;
\\
\textbf{(ETD4RK)}~four-step scheme:
\begin{subequations}
\label{eq:ETD04ad}
\begin{align}
\mathbf{a}&=e^{\mathbf{L}\tau/2}\cdot\mathbf{u}(t) +\mathbf{L}^{-1}\cdot(e^{\mathbf{L}\tau/2}-\mathbf{I})\cdot\mathbf{f}\big(\mathbf{u}(t),t\big)\,,
\label{eq:ETD07}
\\
\mathbf{b}&=e^{\mathbf{L}\tau/2}\cdot\mathbf{u}(t) +\mathbf{L}^{-1}\cdot(e^{\mathbf{L}\tau/2}-\mathbf{I})\cdot\mathbf{f}\big(\mathbf{a},t+\tau/2\big)\,,
\label{eq:ETD08}
\\
\mathbf{c}&=e^{\mathbf{L}\tau/2}\cdot\mathbf{a} +\mathbf{L}^{-1}\cdot(e^{\mathbf{L}\tau/2}-\mathbf{I})\cdot \left[2\mathbf{f}\big(\mathbf{b},t+\tau/2\big)-\mathbf{f}\big(\mathbf{u}(t),t\big)\right],
\label{eq:ETD09}
\\
\mathbf{u}(t+\tau)&=e^{\mathbf{L}\tau}\cdot\mathbf{u}(t) +\frac{\mathbf{L}^{-3}}{\tau^2}\cdot\Big(\big[-4-\mathbf{L}\tau
+e^{\mathbf{L}\tau}(4-3\mathbf{L}\tau
+\mathbf{L}^2\tau^2)\big]\cdot\mathbf{f}\big(\mathbf{u}(t),t\big)
\nonumber\\
&\qquad\qquad\qquad\qquad\qquad{}
+2\left[2+\mathbf{L}\tau +e^{\mathbf{L}\tau}(-2+\mathbf{L}\tau)\right]\cdot \big(\mathbf{f}(\mathbf{a},t+\tau/2)+\mathbf{f}(\mathbf{b},t+\tau/2)\big)
\nonumber\\
&\qquad\qquad\qquad\qquad\qquad{}
+\left[-4-3\mathbf{L}\tau-\mathbf{L}^2\tau^2 +e^{\mathbf{L}\tau}(4-\mathbf{L}\tau)\right]\cdot\mathbf{f}\big(\mathbf{c},t+\tau\big)
\Big)\,.
\label{eq:ETD10}
\end{align}
\end{subequations}
The error of the ETD2RK scheme on one time step is $-\tau^3\ddot{\mathbf{f}}/12$~\cite{Cox-Matthews-2002}, the ETD3RK scheme error is $\propto\tau^4\mathrm{d}^3\mathbf{f}/\mathrm{d}t^3$, the ETD4RK scheme error is $\propto\tau^5\mathrm{d}^4\mathbf{f}/\mathrm{d}t^4$.

In the numerical schemes above, the reciprocal matrix $\mathbf{L}^{-1}$ is introduced for the brevity of expressions using the exponential; wherever the matrix $\mathbf{L}^{-n}$ appears in equations, it is multiplied by expressions the series representation of which starts from the matrix $\mathbf{L}^{m}$, $m\ge n$.
For instance, $\mathbf{L}^{-1}\cdot(e^{\mathbf{L}\tau}-\mathbf{I})$ is actually a formal expression for the series $\tau\mathbf{I}+\tau^2\mathbf{L}/2+\tau^3\mathbf{L}^2/3!+\tau^4\mathbf{L}^3/4!+\dots$\,. Therefore, the zero eigenvalues of the matrix $\mathbf{L}$, which would make $\mathbf{L}^{-1}\to\infty$, are not an issue; the divergence of $\mathbf{L}^{-1}$ only prevents the usage of a shorter form of equations. Albeit, in this paper, we are handling the case of a nondiagonal shape of matrix $\mathbf{L}$, for which analytical calculation of the matrix $e^{\mathbf{L}\tau}$ can be impossible or problematic.

More convenient for our approach is to recast the numerical schemes in terms of the solutions of the following auxiliary problems~\cite{Permyakova-Goldobin-2020}. The general solution to the Cauchy problem
\begin{equation}
\dot{\mathbf{u}}=\mathbf{L}\cdot\mathbf{u}
\label{eq:ETD11}
\end{equation}
is given by the matrix $\mathbf{Q}(\tau)\equiv e^{\mathbf{L}\tau}$ such that
\[
\mathbf{u}\big(t=\tau|\mathbf{u}(0),\mathbf{f}=0\big)=\mathbf{Q}(\tau)\cdot\mathbf{u}(0)\,.
\]
Hence, one obtains the definition of the matrix $\mathbf{Q}(\tau)$, which is convenient for numerical simulations:
\begin{equation}
Q_{jk}(\tau)=u_j\big(t=\tau|u_l(0)=\delta_{lk},\mathbf{f}=0\big)\,.
\label{eq:ETD12}
\end{equation}
The solutions to the problems
\begin{equation}
\dot{\mathbf{u}}=\mathbf{L}\cdot\mathbf{u}+\mathbf{g}\,t^n,\qquad
\mathbf{u}(0)=0\,,\qquad
\mathbf{g}=const
\label{eq:ETD13}
\end{equation}
are given by the matrices $\mathbf{M}_n(\tau)\equiv\int_0^\tau e^{\mathbf{L}(\tau-t)}t^{n-1}\mathrm{d}t$\,:
\[
\mathbf{u}\big(t=\tau|\mathbf{u}(0)=0,\mathbf{f}(t)=\mathbf{g}t^{n-1}\big)=\mathbf{M}_n(\tau)\cdot\mathbf{g}\,.
\]
Hence, one obtains the definition of the matrix $\mathbf{M}_n(\tau)$, which is convenient for numerical simulations:
\begin{equation}
\big(\mathbf{M}_n(\tau)\big)_{jk}=u_j\big(t=\tau|\mathbf{u}(0)=0,f_l(t)=\delta_{lk}t^{n-1}\big)\,.
\label{eq:ETD14}
\end{equation}
Using the recurrent relationship $\mathbf{M}_{n+1}(\tau)=\mathbf{L}^{-1}\cdot(-\tau^n\mathbf{I}+n\mathbf{M}_n(\tau))$, one can obtain
\begin{subequations}
\label{eq:ETDM13}
\begin{align}
\mathbf{M}_1(\tau)&=\mathbf{L}^{-1}\cdot(e^{\mathbf{L}\tau}-\mathbf{I})\,,
\label{eq:ETDM1}\\
\mathbf{M}_2(\tau)&=\mathbf{L}^{-2}\cdot(e^{\mathbf{L}\tau}-\mathbf{I}-\mathbf{L}\tau)\,,
\label{eq:ETDM2}\\
\mathbf{M}_3(\tau)&=\mathbf{L}^{-3}\cdot(2e^{\mathbf{L}\tau}-2\mathbf{I}-2\mathbf{L}\tau-\mathbf{L}^2\tau^2)\,,
\label{eq:ETDM3}\\
&\dots\;.
\nonumber
\end{align}
\end{subequations}

As proposed in~\cite{Permyakova-Goldobin-2020}, for an essentially nondiagonal shape of matrix $\mathbf{L}$, the evaluation of matrices $\mathbf{Q}$ and $\mathbf{M}_n$ defined by equations~(\ref{eq:ETD12}) and (\ref{eq:ETD14}) can be conducted via the direct numerical integration of problems~(\ref{eq:ETD11}) and (\ref{eq:ETD13}) with a very small time stepsize but on a short time interval --- one step of the ETD scheme $\tau$.

In terms of matrices $\mathbf{M}_n$ the schemes (\ref{eq:ETD02ab}), (\ref{eq:ETD03ac}), and (\ref{eq:ETD04ad}) acquire the following form:
\\
\textbf{(ETD2RK)}:
\begin{subequations}
\label{eq:ETD09ab}
\begin{align}
\mathbf{a}&=\mathbf{Q}\cdot\mathbf{u}(t) +\mathbf{M}_1\cdot\mathbf{f}\big(\mathbf{u}(t),t\big)\,,
\label{eq:ETD15}
\\
\mathbf{u}(t+\tau)&=\mathbf{a} +\tau^{-1}\mathbf{M}_2\cdot\left[\mathbf{f}\big(\mathbf{a},t+\tau\big) -\mathbf{f}\big(\mathbf{u}(t),t\big)\right]\,;
\label{eq:ETD16}
\end{align}
\end{subequations}
\\
\textbf{(ETD3RK)}:
\begin{subequations}
\label{eq:ETD10ac}
\begin{align}
\mathbf{a}&=\mathbf{Q}_\frac{\tau}{2}\cdot\mathbf{u}(t) +\mathbf{M}_{1,\frac{\tau}{2}}\cdot\mathbf{f}\big(\mathbf{u}(t),t\big)\,,
\label{eq:ETD17}
\\
\mathbf{b}&=\mathbf{Q}\cdot\mathbf{u}(t) +\mathbf{M}_1\cdot \left[2\mathbf{f}\big(\mathbf{a},t+\tau/2\big)-\mathbf{f}\big(\mathbf{u}(t),t\big)\right]\,,
\label{eq:ETD18}
\\
\mathbf{u}(t+\tau)&=\mathbf{Q}\cdot\mathbf{u}(t) +\left[\frac{2\mathbf{M}_3}{\tau^2}-\frac{3\mathbf{M}_2}{\tau}+\mathbf{M}_1\right]\cdot\mathbf{f}\big(\mathbf{u}(t),t\big)
\nonumber\\
&\qquad{}
-\left[\frac{4\mathbf{M}_3}{\tau^2}-\frac{4\mathbf{M}_2}{\tau}\right]\cdot\mathbf{f}\big(\mathbf{a},t+\tau/2\big)
+\left[\frac{2\mathbf{M}_3}{\tau^2}-\frac{\mathbf{M}_2}{\tau}\right]\cdot\mathbf{f}\big(\mathbf{b},t+\tau\big)\,,
\label{eq:ETD19}
\end{align}
\end{subequations}
where $\mathbf{Q}_\frac{\tau}{2}=\mathbf{Q}(\tau/2)$, $\mathbf{M}_{1,\frac{\tau}{2}}=\mathbf{M}_1(\tau/2)$, and for the full stepsize $\tau$ the corresponding subscript is omitted;
\\
\textbf{(ETD4RK)}:
\begin{subequations}
\label{eq:ETD11ad}
\begin{align}
\mathbf{a}&=\mathbf{Q}_\frac{\tau}{2}\cdot\mathbf{u}(t) +\mathbf{M}_{1,\frac{\tau}{2}}\cdot\mathbf{f}\big(\mathbf{u}(t),t\big)\,,
\label{eq:ETD20}
\\
\mathbf{b}&=\mathbf{Q}_\frac{\tau}{2}\cdot\mathbf{u}(t) +\mathbf{M}_{1,\frac{\tau}{2}}\cdot\mathbf{f}\big(\mathbf{a},t+\tau/2\big)\,,
\label{eq:ETD21}
\\
\mathbf{c}&=\mathbf{Q}_\frac{\tau}{2}\cdot\mathbf{a} +\mathbf{M}_{1,\frac{\tau}{2}}\cdot \left[2\mathbf{f}\big(\mathbf{b},t+\tau/2\big)-\mathbf{f}\big(\mathbf{u}(t),t\big)\right]\,,
\label{eq:ETD22}
\\
\mathbf{u}(t+\tau)&=\mathbf{Q}\cdot\mathbf{u}(t) +\left[\frac{2\mathbf{M}_3}{\tau^2}-\frac{3\mathbf{M}_2}{\tau}+\mathbf{M}_1\right]\cdot\mathbf{f}\big(\mathbf{u}(t),t\big)
\nonumber\\
&
-\left[\frac{2\mathbf{M}_3}{\tau^2}-\frac{2\mathbf{M}_2}{\tau}\right]\cdot\big(\mathbf{f}(\mathbf{a},t+\tau/2)+\mathbf{f}(\mathbf{b},t+\tau/2)\big)
+\left[\frac{2\mathbf{M}_3}{\tau^2}-\frac{\mathbf{M}_2}{\tau}\right]\cdot\mathbf{f}\big(\mathbf{c},t+\tau\big)\,.
\label{eq:ETD23}
\end{align}
\end{subequations}

In the literature, for the difficult case of non-diagonalizable $\mathbf{L}$, the matrices $\mathbf{Q}$ and $\mathbf{M}_n$ are calculated approximately from formulas~(\ref{eq:ETDM13}) either by means of the Taylor expansions and projection into Krylov subspaces~\cite{Hochbruck-Lubich-Selhofer-1998,Tokman-2006} or with discretized contour integrals on the complex plane~\cite{Kassam-Trefethen-2005}. Additionally, for small eigenvalues of matrix $\tau\mathbf{L}$ (the presence of which is abundant), the computations with formula $(e^{\tau\mathbf{L}}-\mathbf{I})/(\tau\mathbf{L})$ (\ref{eq:ETDM1}) suffer from the cancellation errors which rapidly grow for $\mathbf{M}_n$ as $n$ increases further above $1$. Tackling the latter issue is involved but doable~\cite{Cox-Matthews-2002,Kassam-Trefethen-2005}; this issue does not appear with the algorithm we consider in this paper.
Moreover, the technical implementation of this algorithm and the text below do not employ any advanced topics from Linear algebra and Tensor calculus.

\begin{figure*}[t]
\centerline{
\includegraphics[width=0.47\textwidth]{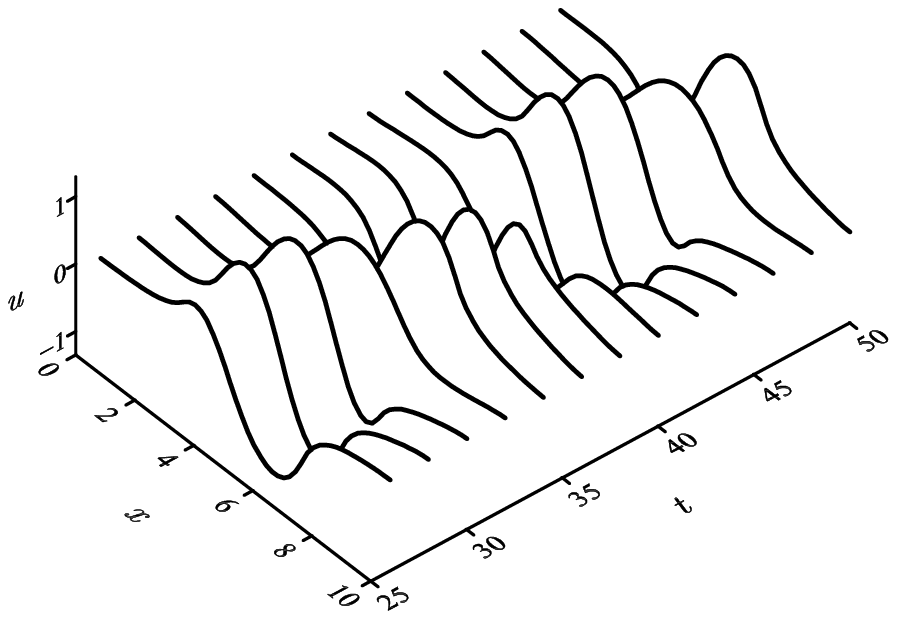}
\qquad
\includegraphics[width=0.47\textwidth]{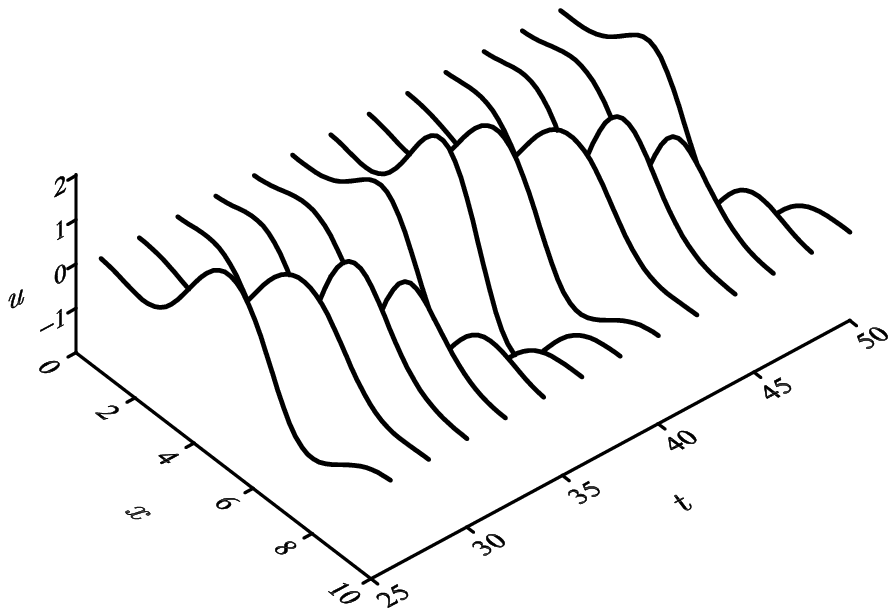}
}
\caption{\label{fig1} Oscillatory solutions of Cahn--Hilliard equation~(\ref{eq:CH01}) (left panel) and Matthews--Cox equation~(\ref{eq:MC01}) (right panel) are plotted for advection velocity $v=1$ and localized excitability $q(x)$ (\ref{eq:CH02}) in the domain of length $L=10$.}
\end{figure*}

\section{Examples \label{sec:Examples}}
\subsection{1D Cahn--Hilliard equation \label{sec:CHE}}

For a broad class of active media systems the pattern formation in thin layers is governed by the Cahn--Hilliard equation~\cite{Knobloch-1990,Shtilman-Sivashinsky-1991,Goldobin-Shklyaeva-2008} (CHE), which as well describes the spinodal decomposition of two-component mixtures~\cite{Cahn-Hilliard-1958,Golovin-etal-2001,Kuramoto-Tsuzuki-1976}. For the generality of consideration, we will also allow for an advective transfer along the layer~\cite{Golovin-etal-2001,Goldobin-Shklyaeva-2008}. In the one-dimensional case (the patterns are homogeneous along the second direction in the layer plane), CHE for the field $u(x,t)$ with advection velocity $v$ reads
\begin{equation}
\frac{\partial u}{\partial t}=-v\frac{\partial u}{\partial x}
-\frac{\partial^2}{\partial x^2}\left[q(x)u+\frac{\partial^2u}{\partial x^2} -u^3\right]\,,
\label{eq:CH01}
\end{equation}
where $q(x)$ is the local deviation of the bifurcation parameter from the instability threshold of the homogeneous infinite layer (the patterns are excited for $q>0$).

For an example, we consider CHE in a domain $0<x<L$ with trivial boundary conditions:
\[
u(0)=\left.\frac{\partial u}{\partial x}\right|_{x=0} =u(L)=\left.\frac{\partial u}{\partial x}\right|_{x=L}=0\,.
\]
The local excitability parameter $q(x)$ is assumed to be positive within certain excitation zone and negative beyond it:
\begin{equation}
q(x)=\left\{
\begin{array}{cc}
2.5\,, & \mbox{  for  }\;\frac{3L}{10}<x<\frac{7L}{10}\,; \\[5pt]
-3\,, & \mbox{  otherwise}\,.
\end{array}
\right.
\label{eq:CH02}
\end{equation}
For a localized excitation $q(x)$, strong enough advection $v$ results in an oscillatory behavior (e.g., see~\cite{Goldobin-2019}). For $v=1$ and $q(x)$ given by~(\ref{eq:CH02}), the oscillatory solution of Eq.~(\ref{eq:CH01}) is presented in Fig.~\ref{fig1}.

\begin{figure*}[t]
\centerline{
\includegraphics[width=0.35\textwidth]{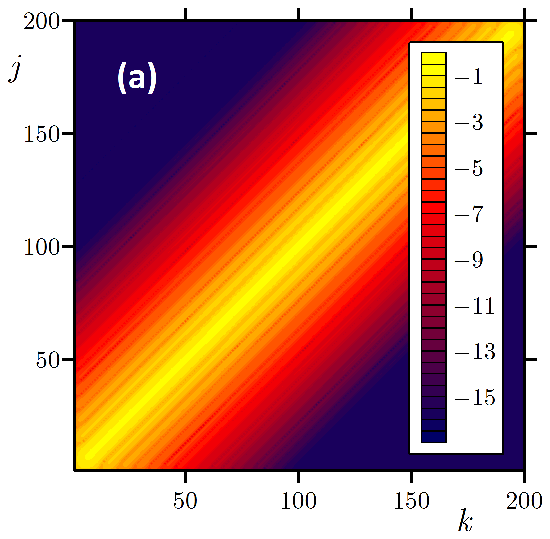}
\qquad\qquad
\includegraphics[width=0.35\textwidth]{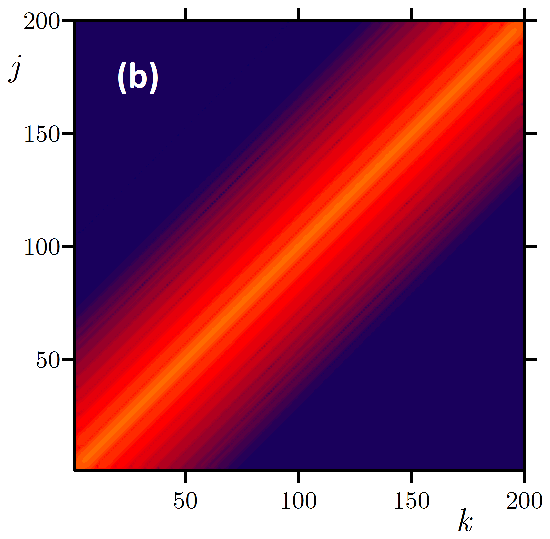}
}
\vspace{7pt}
\centerline{
\includegraphics[width=0.35\textwidth]{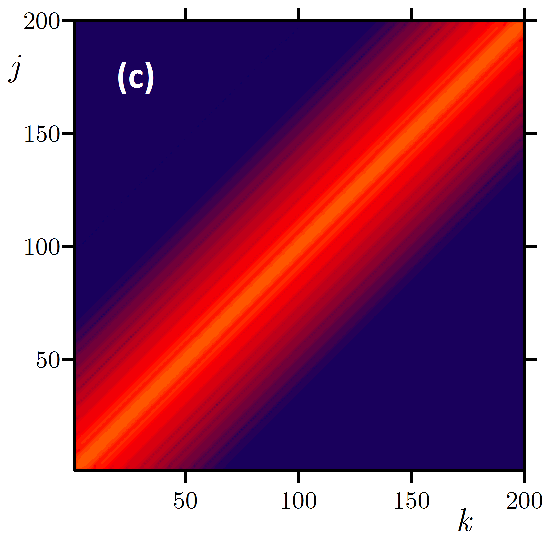}
\qquad\qquad
\includegraphics[width=0.35\textwidth]{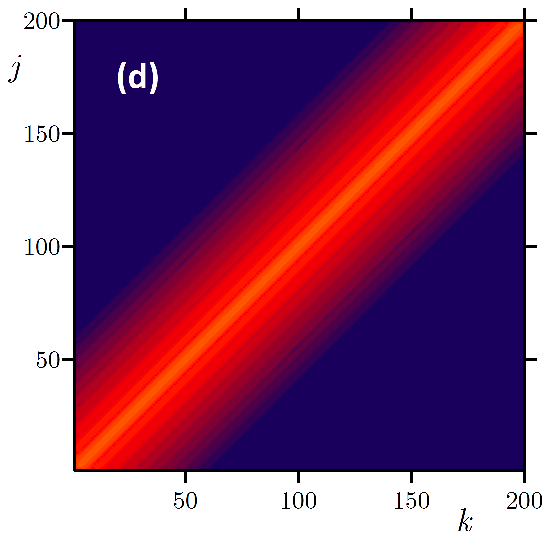}}
\caption{\label{fig2} For Cahn--Hilliard equation~(\ref{eq:CH01}), matrices
$\log_{10}|Q_{jk}|$ (panel \textit{a}),
$\log_{10}|(\mathbf{M}_1)_{jk}|$ (panel \textit{b}),
$\log_{10}|(\tau^{-1}\mathbf{M}_2)_{jk}|$ (panel \textit{c}),
$\log_{10}|(\tau^{-2}\mathbf{M}_3)_{jk}|$ (panel \textit{d})
are plotted versus indices $k$ and $j$ for $\tau=2.5\times10^{-4}$, $L=10$, $N=200$;
$\tau_1=0.1h_x^4=6.25\times10^{-7}$.
}
\end{figure*}

\subsubsection{\label{sec:CH_num_sch}Numerical scheme}

The spatially-discrete version of CHE~(\ref{eq:CH01}) for vector
$\mathbf{u}(t)\equiv\{u_j(t)|j=1,2,...,N\}$, where $u_j(t)=u(jh_x,t)$, $h_x=L/N$, can be written as Eq.~(\ref{eq:ETD01}) with
\begin{subequations}
\label{eq:CH03ab}
\begin{align}
(\mathbf{L}\cdot\mathbf{u})_j&=v\frac{u_{j-1}-u_{j+1}}{2h_x} -\frac{q_{j+1}u_{j+1}-2q_ju_j+q_{j-1}u_{j-1}}{h_x^2}
-\frac{u_{j+2}-4u_{j+1}+6u_j-4u_{j-1}+u_{j-2}}{h_x^4}\;,
\label{eq:CH03}
\end{align}
\begin{equation}
f_j(\mathbf{u},t)=\frac{u_{j+1}^3-2u_j^3+u_{j-1}^3}{h_x^2}\;,
\label{eq:CH04}
\end{equation}
\end{subequations}
where $q_j=q(jh_x)$ (\ref{eq:CH02}).
The trivial boundary conditions we consider imply that one can formally set $u_{-1}=u_0=u_{N+1}=u_{N+2}=0$ in the right-hand parts of (\ref{eq:CH03ab}).

The problem (\ref{eq:ETD01},\ref{eq:CH03ab}) can be numerically simulated with the basic Euler scheme or the predictor--corrector (PC) one. For both of these schemes, the forth-order $x$-derivative term imposes the limitation on the time stepsize, $\tau_1<h_x^4/8$; for a bigger time stepsize the direct numerical simulation becomes unstable. More sophisticated Runge--Kutta type schemes or any other high-order methods are not needed here as the PC scheme with such a small time stepsize already warrants an excessive accuracy of the numerical integration in time. Henceforth, we will employ the PC scheme:
\begin{subequations}
\label{eq:CH05ad}
\begin{align}
\mathbf{F}(t)&=\mathbf{L}\cdot\mathbf{u}(t)+\mathbf{f}\big(\mathbf{u}(t),t\big)\;,
\label{eq:CH05}
\\
\mathbf{a}&=\mathbf{u}(t)+\mathbf{F}(t)\tau_1\;,
\label{eq:CH06}
\\
\mathbf{F}(t+\tau_1)&=\mathbf{L}\cdot\mathbf{a}+\mathbf{f}\big(\mathbf{a},t\big)\;,
\label{eq:CH07}
\\
\mathbf{u}(t+\tau_1)&=\mathbf{u}(t)+\frac{\mathbf{F}(t)+\mathbf{F}(t+\tau_1)}{2}\tau_1\;.
\label{eq:CH08}
\end{align}
\end{subequations}

For problem (\ref{eq:ETD01},\ref{eq:CH03ab}), the PC scheme~(\ref{eq:CH05ad}) was employed to integrate the auxiliary problems (\ref{eq:ETD11}) and (\ref{eq:ETD13}) for $t\in[0,\tau]$ to evaluate matrices $\mathbf{Q}_\frac{\tau}{2}$, $\mathbf{Q}$ (\ref{eq:ETD12}) and $\mathbf{M}_{1,\frac{\tau}{2}}$, $\mathbf{M}_n$, $n=1,2,3$ (\ref{eq:ETD14}), respectively. The time stepsize for the direct numerical integration of the auxiliary problems was $\tau_1=0.1h_x^4$. In Fig.~\ref{fig2}, a sample structure of these matrices can be seen.

\begin{figure*}[t]
\centerline{
\includegraphics[width=0.39\textwidth]{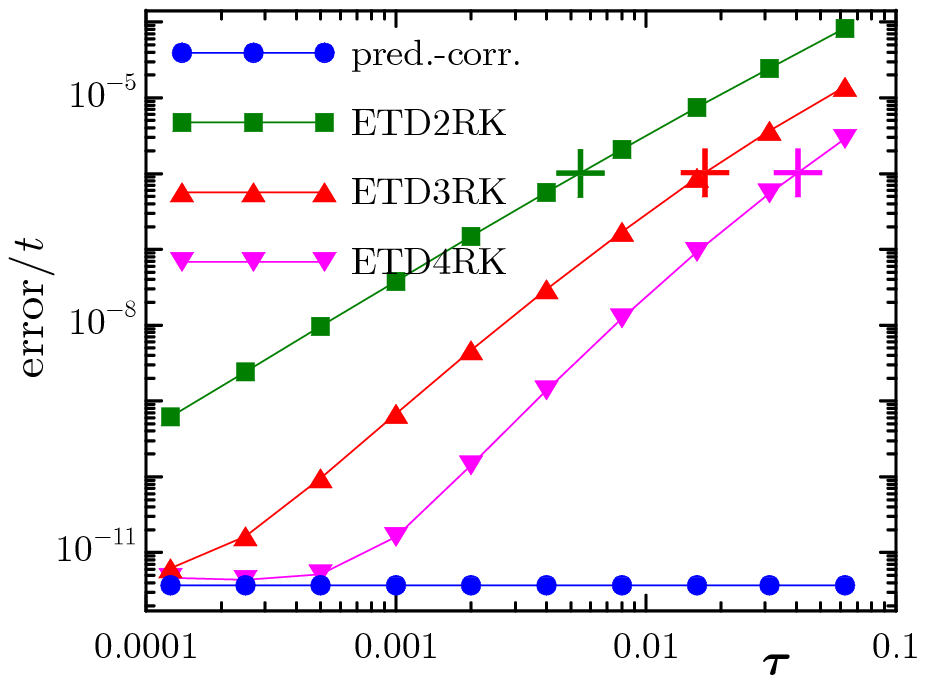}
\qquad\qquad\qquad
\includegraphics[width=0.438\textwidth]{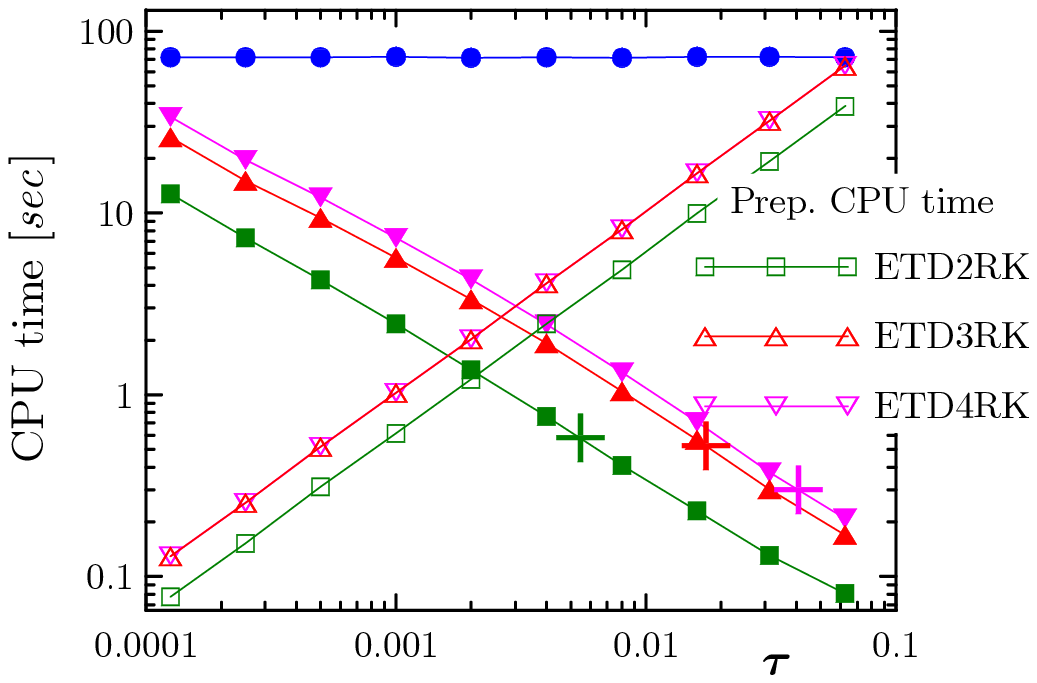}}
\caption{\label{fig3} The error rate and the CPU time for the numerical simulation of Cahn--Hilliard equation~(\ref{eq:CH01}) with several ETD schemes over the time interval $0<t<50$ are plotted versus the ETD stepsize $\tau$ for $N=200$ (i.e.\ $h_x=0.05$).
Solid squares: ETD2RK, solid up-triangles: ETD3RK, solid down-triangles: ETD4RK, solid circles: predictor--corrector scheme with fixed time stepsize $\tau_1=0.1h_x^4=6.25\times10^{-7}$.
CPU times for the preliminary  (preparatory) calculation of the matrices $\mathbf{Q}$ and $\mathbf{M}_n$ required for the respective ETD scheme are plotted in the right panel with open symbols. CPU times are provided for the processor
Intel(R) Core(TM) i7-4790K CPU 4.00 GHz, disabled hypertrading; RAM: DDR3 16 GB; program in FORTRAN.}
\end{figure*}
\begin{figure*}[t]
\centerline{
\includegraphics[width=0.39\textwidth]{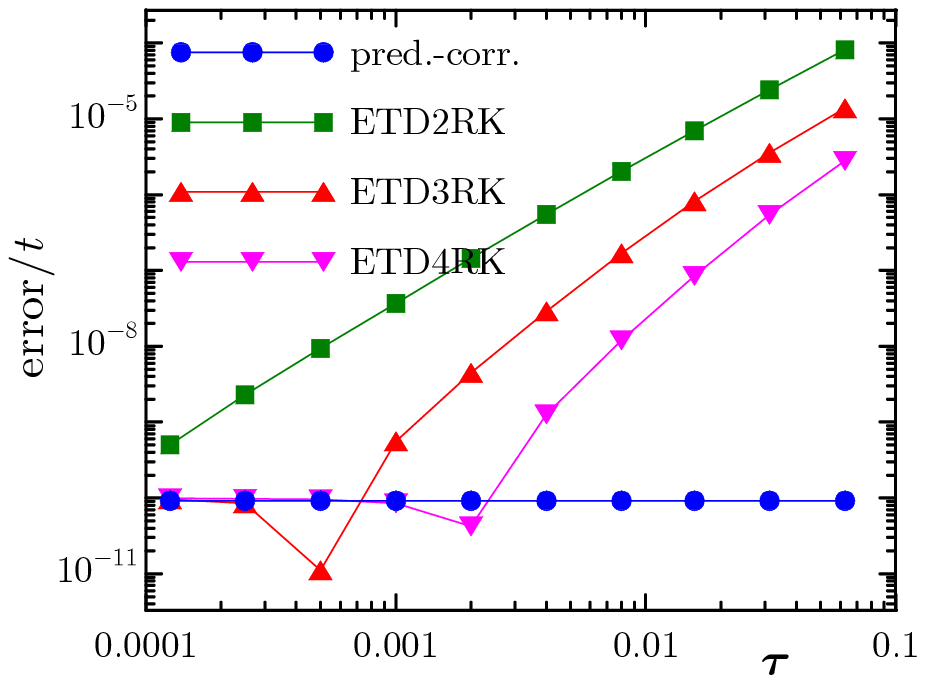}
\qquad\qquad\qquad
\includegraphics[width=0.438\textwidth]{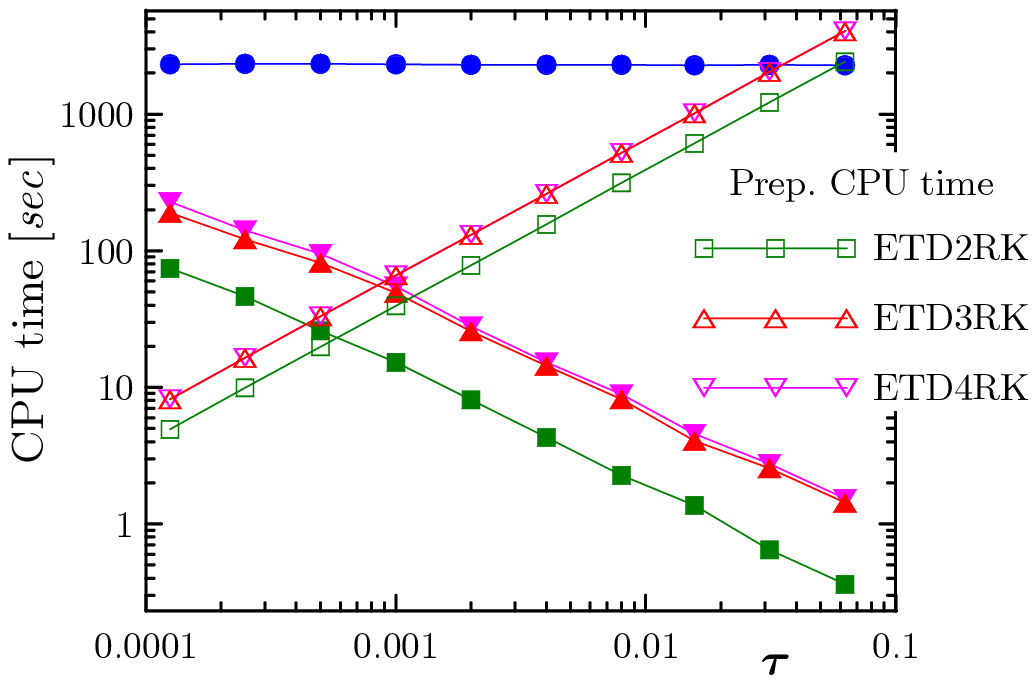}}
\caption{\label{fig4} The error rate and the CPU time for the numerical simulation of CHE~(\ref{eq:CH01}) with several ETD schemes are plotted versus the ETD stepsize $\tau$ for $N=400$ (i.e.\ $h_x=0.025$);
$\tau_1=0.1h_x^4=3.90625\times10^{-8}$.
See Caption to Fig.~\ref{fig3} for notations and values of other parameters.}
\end{figure*}

\subsubsection{Accuracy and performance of numerical simulation with exponential time differencing schemes}

The ETD2RK, ETD3RK, and ETD4RK schemes with matrices $\mathbf{Q}_\frac{\tau}{2}$, $\mathbf{Q}$, $\mathbf{M}_{1,\frac{\tau}{2}}$, $\mathbf{M}_n$, $n=1,2,3$ calculated as described in Sec.~\ref{sec:CH_num_sch} were employed for the numerical simulation of CHE~(\ref{eq:CH01}) with $v=1$ and $q(x)$ given by (\ref{eq:CH02}) in the domain of length $L=10$  (see Fig.~\ref{fig1}). In Figs.~\ref{fig3} and \ref{fig4}, the dependencies of the accuracy and the simulation performance versus the stepsize $\tau$ of the ETD schemes are presented.

In Fig.~\ref{fig3}, one can compare the performance of ETD methods to that of the PC method (blue solid circles). The numerical error of a specific scheme is evaluated as the deviation of the result from the result calculated with the PC scheme~(\ref{eq:CH05ad}) with $\tau_1=0.01h_x^4$, which is by a factor of 10 smaller then the time stepsize $0.1h_x^4$ needed for practical computations. The accuracy of the PC method is excessive, but it is dictated by the maximal admissible value of the time stepsize $h_x^4/8$. In turn, the ETD methods can provide decent accuracy even for large values of $\tau$, where they perform much faster (see the right panel) than the PC method. For example, with a required accuracy of $10^{-6}$ (see crests in Fig.~\ref{fig3}), one can pick-up as large time stepsize as $\tau\approx0.005$ for the ETD2RK scheme, $0.02$ for ETD3RK, and $0.04$ for ETD4RK. The biggest performance gain is achieved with ETD4RK scheme and exceeds a factor of $100$ for all three schemes. For an increased accuracy of the spatial discretization (see Fig.~\ref{fig4} where $N$ is doubled compared to the case of Fig.~\ref{fig3}), the performance gain provided by ETD schemes significantly increases.

Note, the CPU time for the preliminary calculations of $\mathbf{Q}$ and $\mathbf{M}_n$ increases linearly with $\tau$ and can become nonsmall (see open symbols in the right panels of Figs.~\ref{fig3} and \ref{fig4}). This CPU time can be significant for the problems where simulation over short or moderate intervals of time $t$ is sufficient. For these problems, an optimal performance is achieved for the values of $\tau$, which provide approximately equal CPU times for the preliminary calculation of $\mathbf{Q}$ and $\mathbf{M}_n$ and the run of the ETD scheme (see the crossings of the lines marked by solid and open symbols in Figs.~\ref{fig3} and \ref{fig4}). The overall performance gain is decreased, but still large. The ETD2RK scheme becomes the most efficient one in this case. However, the most computationally demanding problems are those of complex dynamics~\cite{Golovin-etal-2001,Watson-etal-2003,Podolny-etal-2005,Shklyaev-etal-2012,Samoilova-Nepomnyashchy-2019,Samoilova-Nepomnyashchy-2020,Zincenko-Petrovskii-etal-2021,Pal-Petrovskii-etal-2021,Chowdhury-Banerjee-Petrovskii-2022,Boaretto-etal-Kurths-2018,Godavarthi-etal-Kurths-2020,Sakaguchi-Shinomoto-Kuramoto-1988,Klinshov-etal-2021,Franovic-etal-2022} (including spatiotemporal chaos) and dynamics in the presence of frozen parametric disorder~\cite{Hammele-Schuler-Zimmermann-2006,Goldobin-Shklyaeva-2013,Goldobin-Shklyaeva-2009,Goldobin-2010,Goldobin-2019}, where the CPU time for the preliminary calculations can be neglected.

In~\ref{sec:app}, our analysis is complemented with the examination of the accuracy of approximate calculation of $\mathbf{Q}$ and $\mathbf{M}_n$.
The inaccuracies in $\mathbf{Q}$ and $\mathbf{M}_n$ are primarily contributed by fast decaying modes, for which the contribution into the inaccuracy of a numerically simulated solution $\mathbf{u}(t)$ is exponentially suppressed.
The impact of inaccuracies of $\mathbf{Q}$ and $\mathbf{M}_n$ related to different decaying perturbation modes onto the inaccuracy of $\mathbf{u}(t)$ is practically impossible to unpack.
The net error of the numerical method is contributed by two sources: the error of the ETD scheme and the error of approximate calculation of matrices $\mathbf{Q}$ and $\mathbf{M}_n$. With all the importance of the information on the inaccuracy of $\mathbf{Q}$ and $\mathbf{M}_n$, the net error decides the practical applicability of the approach; this error is reported in the left panels of Figs.~\ref{fig3} and \ref{fig4}.

\begin{figure*}[t]
\centerline{
\includegraphics[width=0.35\textwidth]{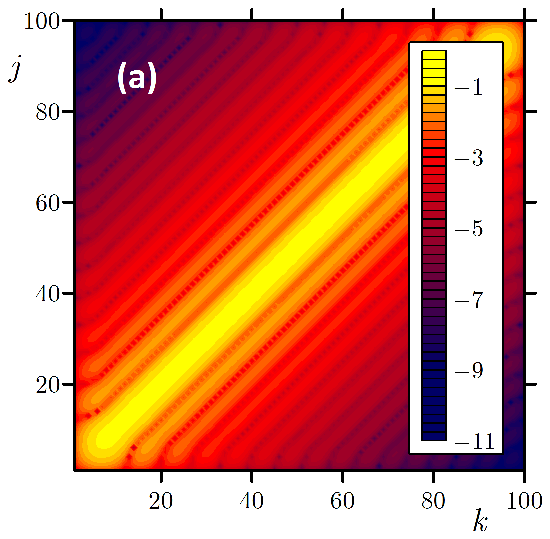}
\qquad\qquad
\includegraphics[width=0.35\textwidth]{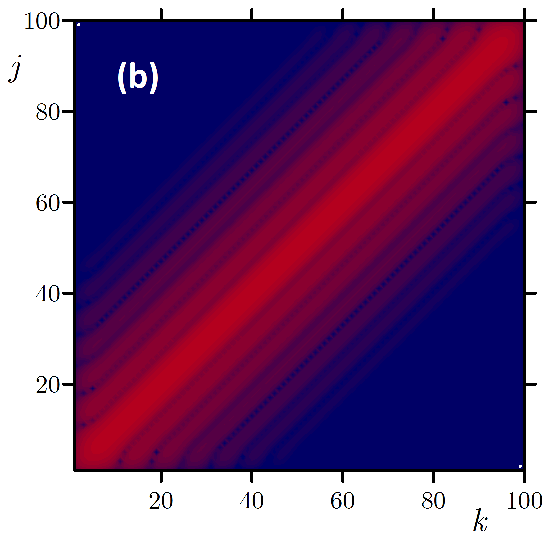}
}
\vspace{7pt}
\centerline{
\includegraphics[width=0.35\textwidth]{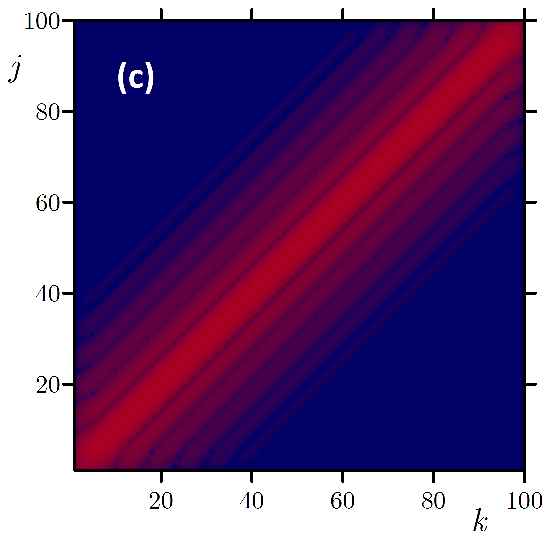}
\qquad\qquad
\includegraphics[width=0.35\textwidth]{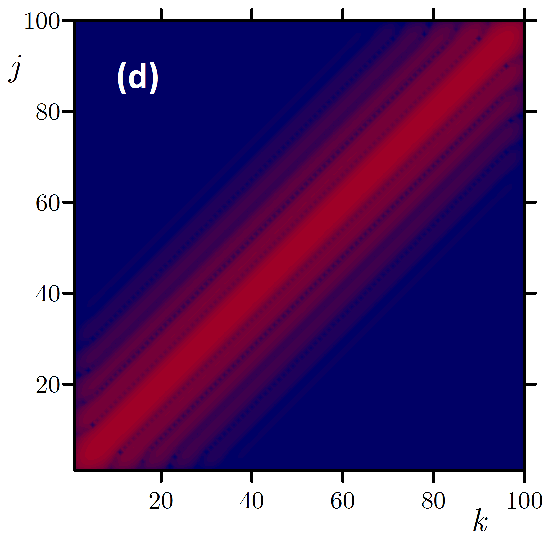}}
\caption{\label{fig5} For Matthews--Cox equation~(\ref{eq:MC01}), matrices
$\log_{10}|Q_{jk}|$ (panel \textit{a}),
$\log_{10}|(\mathbf{M}_1)_{jk}|$ (panel \textit{b}),
$\log_{10}|(\tau^{-1}\mathbf{M}_2)_{jk}|$ (panel \textit{c}),
$\log_{10}|(\tau^{-2}\mathbf{M}_3)_{jk}|$ (panel \textit{d})
are plotted versus indices $k$ and $j$ for $\tau=1.6\times10^{-4}$, $L=10$, $N=100$;
$\tau_1=0.02h_x^6=2\times10^{-8}$.
}
\end{figure*}

\subsection{1D pattern formation with a conservation law \label{sec:MCE}}

The conservation laws (of the chemical specie mass, etc.) are structurally stable (persistent) properties of systems and influence the general form of equations governing the pattern formation. Matthews and Cox~\cite{Matthews-Cox-2000} studied the model equation governing pattern formation with an additional conservation law:
\begin{equation}
\frac{\partial u}{\partial t}=-v\frac{\partial u}{\partial x}
-\frac{\partial^2}{\partial x^2}\left[q(x)u-2\frac{\partial^2u}{\partial x^2} -\frac{\partial^4u}{\partial x^4} -u^3\right]\,,
\label{eq:MC01}
\end{equation}
where the terms inside the brackets are just those of the well-studied Swift--Hohenberg
equation~\cite{Swift-Hohenberg-1977} and an additional advective $v$-term is introduced~\cite{Matthews-Cox-2000b}. Let us consider the Matthews--Cox equation (MCE) as an example of the pattern formation equation with the sixth-order spatial derivative.

We consider MCE in a domain $0<x<L$ with trivial boundary conditions:
\[
u(0)=\left.\frac{\partial u}{\partial x}\right|_{x=0}\!\!\!
=\left.\frac{\partial^2u}{\partial x^2}\right|_{x=0}\!\!\!
=u(L)=\left.\frac{\partial u}{\partial x}\right|_{x=L}\!\!\!
=\left.\frac{\partial^2u}{\partial x^2}\right|_{x=L}\!\!\!
=0\,.
\]
The local excitability parameter $q(x)$ is assumed to be the same as for CHE (\ref{eq:CH02}).
For a localized excitation $q(x)$, strong enough advection $v$ results in an oscillatory behavior. For $v=1$ and $q(x)$ given by~(\ref{eq:CH02}), the oscillatory solution of Eq.~(\ref{eq:MC01}) is presented in Fig.~\ref{fig1}.

\begin{figure*}[t]
\centerline{
\includegraphics[width=0.39\textwidth]{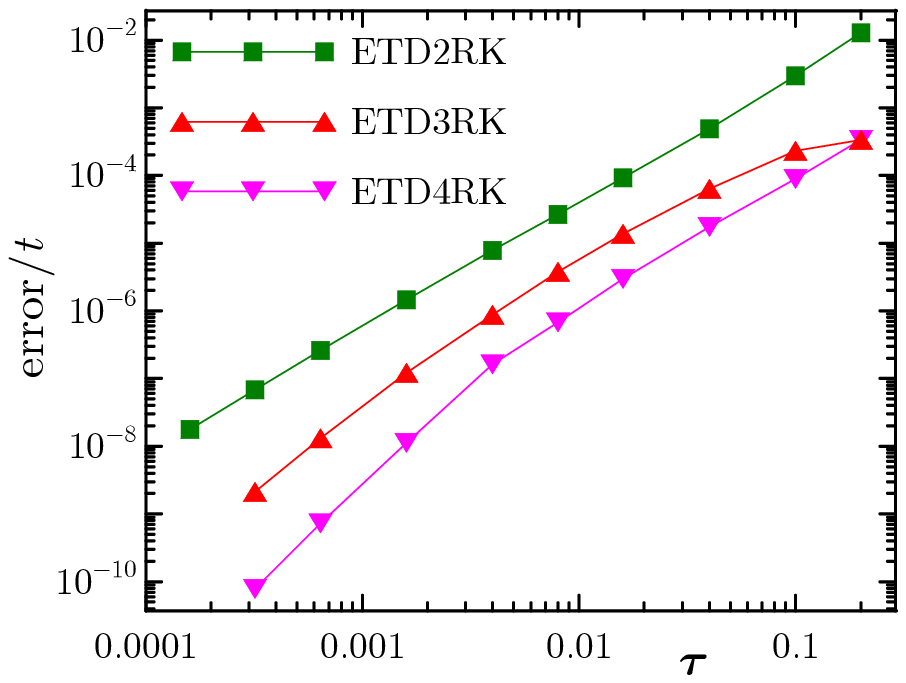}
\qquad\qquad\qquad
\includegraphics[width=0.441\textwidth]{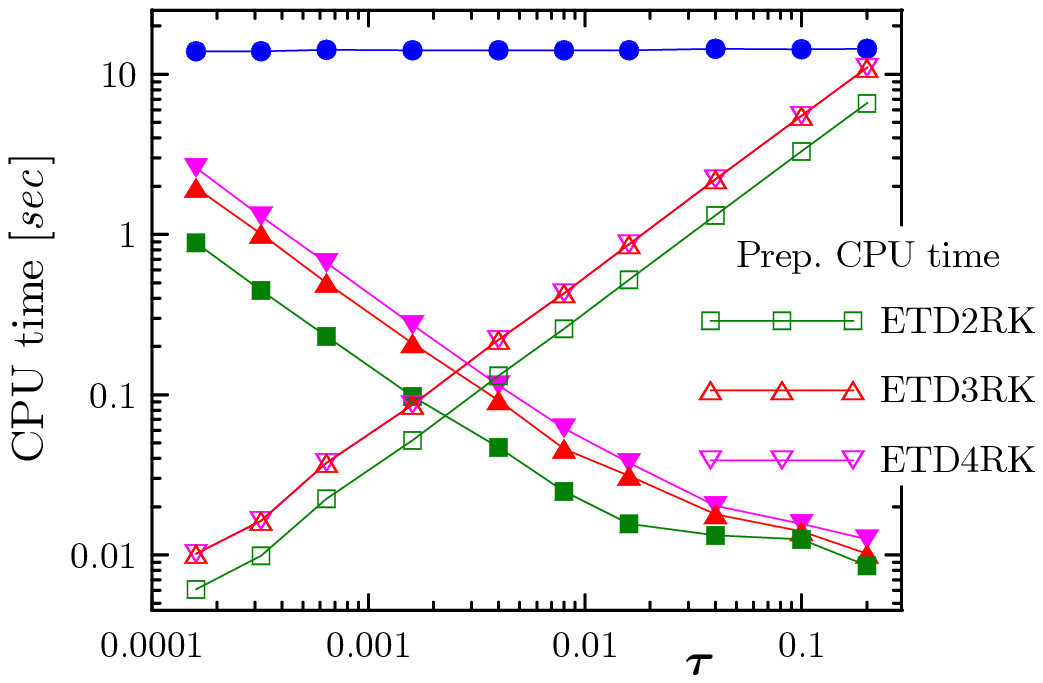}}
\caption{\label{fig6}
The error rate and the CPU time for the numerical simulation of MCE~(\ref{eq:MC01}) with several ETD schemes are plotted versus the ETD stepsize $\tau$ for $N=50$ (i.e.\ $h_x=0.2$);
$\tau_1=0.02h_x^6=1.6\times10^{-7}$.
See Caption to Fig.~\ref{fig3} for notations and values of other parameters.
}
\end{figure*}
\begin{figure*}[t]
\centerline{
\includegraphics[width=0.39\textwidth]{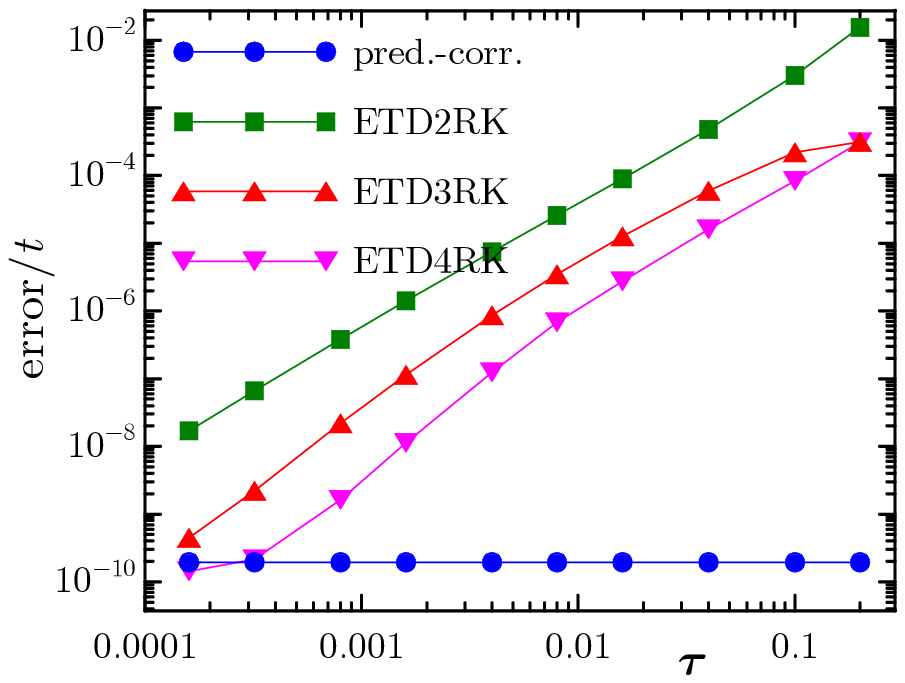}
\qquad\qquad\qquad
\includegraphics[width=0.441\textwidth]{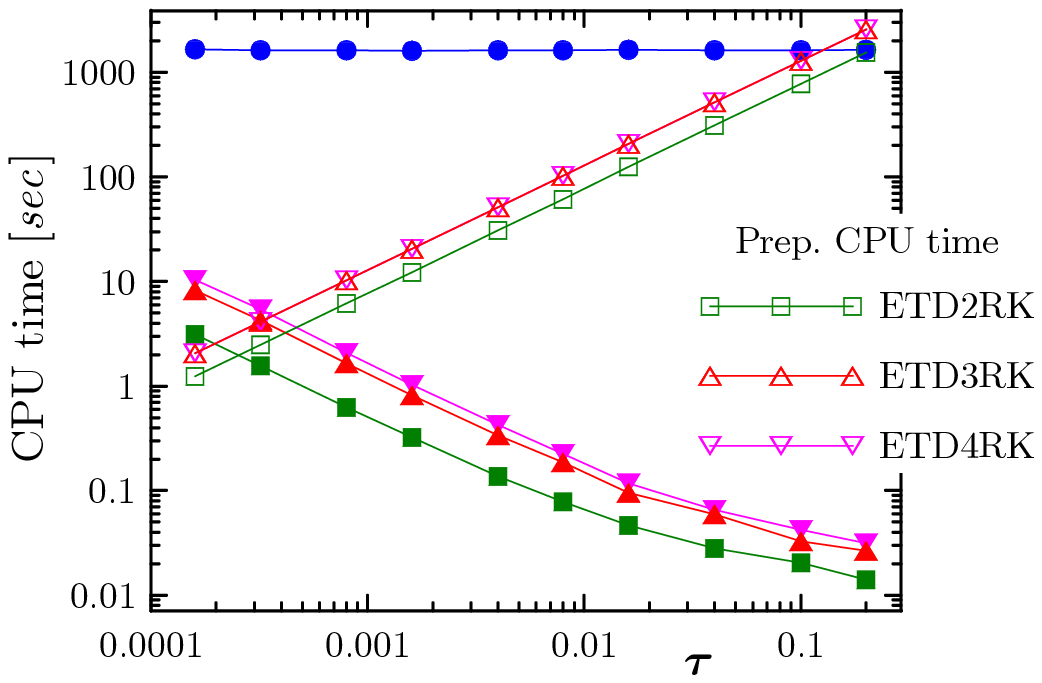}}
\caption{\label{fig7}
The error rate and the CPU time for the numerical simulation of MCE~(\ref{eq:MC01}) with several ETD schemes are plotted versus the ETD stepsize $\tau$ for $N=100$ (i.e.\ $h_x=0.1$);
$\tau_1=0.02h_x^6=2\times10^{-8}$. See Caption to Fig.~\ref{fig3} for notations and values of other parameters.}
\end{figure*}

\subsubsection{\label{sec:MC_num_sch}Numerical scheme}

The spatially-discrete version of MCE~(\ref{eq:MC01}) for vector
$\mathbf{u}(t)\equiv\{u_j(t)|j=1,2,...,N\}$, where $u_j(t)=u(jh_x,t)$, $h_x=L/N$, can be written as Eq.~(\ref{eq:ETD01}) with
\begin{subequations}
\label{eq:MC02ab}
\begin{align}
(\mathbf{L}\cdot\mathbf{u})_j&=v\frac{u_{j-1}-u_{j+1}}{2h_x} -\frac{q_{j+1}u_{j+1}-2q_ju_j+q_{j-1}u_{j-1}}{h_x^2}
+2\frac{u_{j+2}-4u_{j+1}+6u_j-4u_{j-1}+u_{j-2}}{h_x^4}
\nonumber
\\
&\qquad\qquad{}
+\frac{u_{j+3}-6u_{j+2}+15u_{j+1}-20u_j+15u_{j-1}-6u_{j-2}+u_{j-3}}{h_x^6}\;,
\label{eq:MC02a}
\end{align}
\begin{equation}
f_j(\mathbf{u},t)=\frac{u_{j+1}^3-2u_j^3+u_{j-1}^3}{h_x^2}\;,
\label{eq:MC02b}
\end{equation}
\end{subequations}
where $q_j=q(jh_x)$ (\ref{eq:CH02}).
The trivial boundary conditions we consider imply that one can formally set $u_{-2}=u_{-1}=u_0=u_{N+1}=u_{N+2}=u_{N+3}=0$ in the right-hand parts of (\ref{eq:MC02ab}).

The problem (\ref{eq:ETD01},\ref{eq:MC02ab}) can be numerically simulated with the basic Euler scheme or the predictor--corrector one. For both of these schemes, the sixth-order $x$-derivative term imposes the limitation on the time stepsize, $\tau_1<h_x^6/32$; for a bigger time stepsize the direct numerical simulation becomes unstable. Henceforth, we will employ the PC scheme (\ref{eq:CH05ad}).

For problem (\ref{eq:ETD01},\ref{eq:MC02ab}), the PC scheme~(\ref{eq:CH05ad}) was employed to integrate the auxiliary problems (\ref{eq:ETD11}) and (\ref{eq:ETD13}) for $t\in[0,\tau]$ to evaluate matrices $\mathbf{Q}_\frac{\tau}{2}$, $\mathbf{Q}$, $\mathbf{M}_{1,\frac{\tau}{2}}$, $\mathbf{M}_n$, $n=1,2,3$. The time stepsize for the direct numerical integration of the auxiliary problems was $\tau_1=0.02h_x^6$. In Fig.~\ref{fig5}, a sample structure of these matrices can be seen.

\subsubsection{Accuracy and performance of numerical simulation with exponential time differencing schemes}

The ETD2RK, ETD3RK, and ETD4RK schemes with matrices $\mathbf{Q}_\frac{\tau}{2}$, $\mathbf{Q}$, $\mathbf{M}_{1,\frac{\tau}{2}}$, $\mathbf{M}_n$, $n=1,2,3$ calculated as described in Sec.~\ref{sec:MC_num_sch} were employed for the numerical simulation of MCE~(\ref{eq:MC01}) with $v=1$ and $q(x)$ given by (\ref{eq:CH02}) in the domain of length $L=10$ (see Fig.~\ref{fig1}). In Figs.~\ref{fig6} and \ref{fig7}, the dependencies of the accuracy and the simulation performance versus the stepsize $\tau$ of the ETD schemes are presented. One can see even more dramatic performance gain than for the numerical simulation of CHE (\ref{eq:CH01}).

Generally, for the partial differential equations with the highest-order spatial derivative $\partial^m/\partial x^m$, the performance gain for the ETD methods is $\propto h_x^{m-2}$, as anticipated in~\cite{Permyakova-Goldobin-2020}. However, in Figs.~\ref{fig2} and \ref{fig5}, one can see that matrix elements rapidly decay away from the diagonal. For nonlarge $\tau$, the absolute values of a significant fraction of elements are below $\sim10^{-16}$ which is the level of the double precision computer accuracy. The elements smaller than the required error level can be set to zero without any damage to the scheme accuracy. The optimization of a program code for the sparseness of matrices $\mathbf{Q}$ and $\mathbf{M}_n$ allows for a noticeable performance acceleration. For a given stepsize $\tau$, the highest derivative $\partial^m/\partial x^m$ creates a bigger spread of nonsmall values of matrix elements (compare Fig.~\ref{fig2} to Fig.~\ref{fig5}, where the stepsize $\tau$ is even somewhat smaller). Therefore, this additional optimization becomes less beneficial for high-order derivatives, where the gain is actually large without any optimization.

\subsection{Macroscopic dynamics of population of quadratic integrate-and-fire neurons with noise \label{sec:QIFs}}
The next example we consider is the macroscopic dynamics of a large population of coupled oscillators. The quadratic integrate-and-fire neuron (QIF)~\cite{Izhikevich-2007} is not only a useful toy model for the excitable behavior of neurons, but also a normal for the Class~I excitability neurons near the bifurcation point~\cite{Ermentrout-Kopell-1986}. Networks of QIFs can exhibit sophisticated macroscopic dynamics; these dynamics and the effect of endogenic noise on them can be often studied and well understood within the framework of the Fokker--Planck equation~\cite{Ratas-Ryragas-2019,Goldobin-diVolo-Torcini-2021,diVolo-etal-2022,Zheng-Kotani-Jimbo-2021,Goldobin-Permyakova-Klimenko-2024,Pietras-Cestnik-Pikovsky-2023}. However, the problem of reasonably accurate direct numerical simulation of it for some important collective regimes becomes a challenging task. Below we provide the equations for a large recurrent synaptic network of QIFs with global coupling and endogenic Gaussian noise, explain typical accuracy breakdown of its direct numerical simulation, and test three ETD schemes for it.

\subsubsection{Governing equations}
In this section we provide a concise derivation and interpretation for equation system~(\ref{eq:108}) and (\ref{eq:110}), the direct numerical simulation of which can be a challenging task tackled by means of the ETD-methods.

We consider the population of QIFs with endogenic noise:
\begin{align}
&\dot{V}_n=V_n^2+I_n\,,
\label{eq:101}
\\
&I_n=\eta_n+\sigma\zeta_n(t)+Js(t)+I(t)\,,
\label{eq:102}
\end{align}
where $V_n$ is the membrane voltage of the $n$th neuron, $\eta_n$ is the excitability parameter of individual neuron (an isolated neuron is excitable for $\eta_n<0$ and periodically spiking otherwise), $I(t)$ is the external input current, $\sigma\zeta_n(t)$ are independent Gaussian endogenic noises: $\langle\zeta_n(t)\,\zeta_m(t+t^\prime)\rangle=2\delta_{nm}\delta(t^\prime)$, $\delta_{nm}$ is the Kronecker delta which is $1$ for $n=m$ and $0$ otherwise. When $V_n$ reaches $+\infty$ it is reset to $-\infty$ and a synaptic spike is generated~\cite{Ermentrout-Kopell-1986}. The input synaptic current from other neurons $Js(t)$ is characterised by the coupling coefficient $J$, which is negative for an inhibitory coupling, and a common field $s(t)$ equals to the neuron firing rate $r(t)$ for a thermodynamically large population with instantaneous synaptic spikes.

One can introduce a phase variable $\phi$,
\[
V_n=\tan\frac{\phi_n}{2},
\]
and rewrite Eq.~(\ref{eq:101}) as
\begin{equation}
\dot\phi_n=(1-\cos\phi_n)+(1+\cos\phi_n)\left[\eta_n+\sigma\zeta_n(t)+Js(t)+I(t)\right].
\label{eq:103}
\end{equation}
Let us index the QIFs with the value of their parameter $\eta_n$. The Fokker--Planck equation for the probability density $w_\eta(\phi,t)$ of stochastic system~(\ref{eq:103}) reads
\begin{align}
&\frac{\partial w_\eta}{\partial t}+\frac{\partial}{\partial\phi}\Big(\big[1-\cos\phi +(1+\cos\phi)(\eta+Js+I(t))\big]w_\eta\Big)
=\sigma^2\frac{\partial}{\partial\phi}\left((1+\cos\phi)\frac{\partial}{\partial\phi}\big((1+\cos\phi)w_\eta\big)\right).
\label{eq:104}
\end{align}

Now we calculate the firing rate $r$ in terms of $\phi$. Given the distribution of $\eta$ is $g(\eta)$, in the thermodynamic limit of large population, the firing rate equals
 $r(t)=\int q_\eta(\phi=\pi)\,g(\eta)\,\mathrm{d}\eta$,
where $q_\eta$ is the probability density flux,
\begin{align}
q_\eta=\big[1-\cos\phi +(1+\cos\phi)(\eta+Js+I(t))\big]w_\eta
-\sigma^2(1+\cos\phi)\frac{\partial}{\partial\phi}\big((1+\cos\phi)w_\eta\big)\,.
\nonumber
\end{align}
For $\phi=\pi$, the flux $q_\eta(\phi=\pi)=2w_\eta(\pi)$ and one finds
\begin{equation}
r(t)=\int q_\eta(\pi)\,g(\eta)\mathrm{d}\eta=2\int w_\eta(\pi)\,g(\eta)\mathrm{d}\eta\,.
\label{eq:105}
\end{equation}


\begin{figure}[t]
\centerline{
\includegraphics[width=0.46\textwidth]{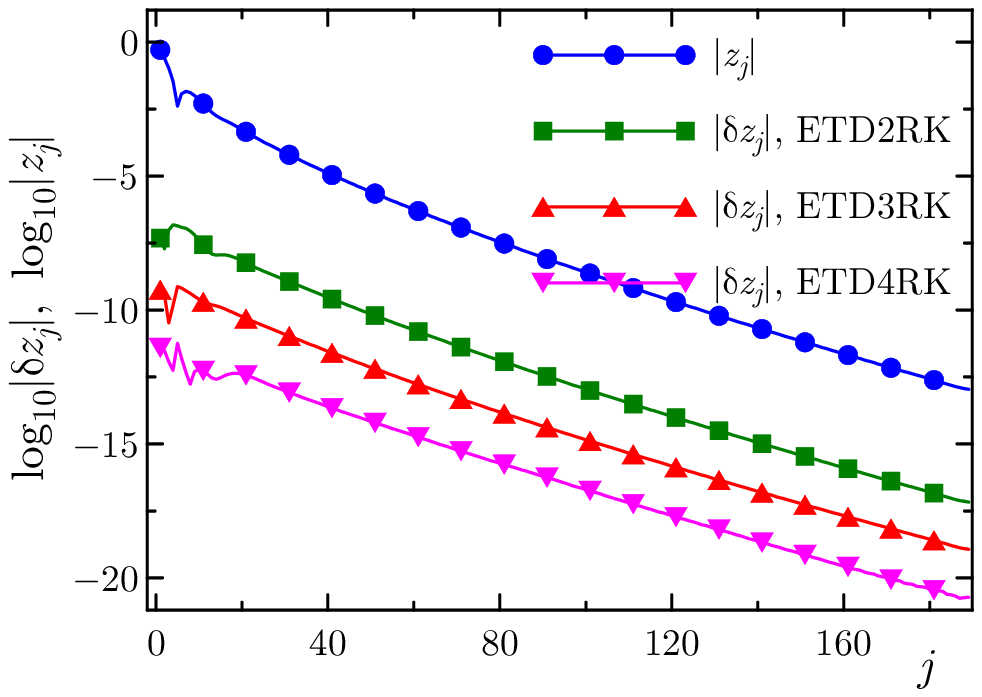}
}
\caption{\label{fig8}
Instantaneous state $\{z_j\}$ of system~(\ref{eq:108}) and (\ref{eq:110}) after a transient process for $\eta_0=-2$, $\gamma=1$, $J=10$, $\sigma^2=10$, and periodically modulated external input current $I(t)=0.3\sin{10t}$. The numerical errors $\{\delta z_j\}$ are plotted for the simulations with ETD schemes with $\tau=0.025$
and $\tau_1=1.25\times10^{-6}$.
}
\end{figure}
\begin{figure*}[t]
\centerline{
\includegraphics[width=0.35\textwidth]{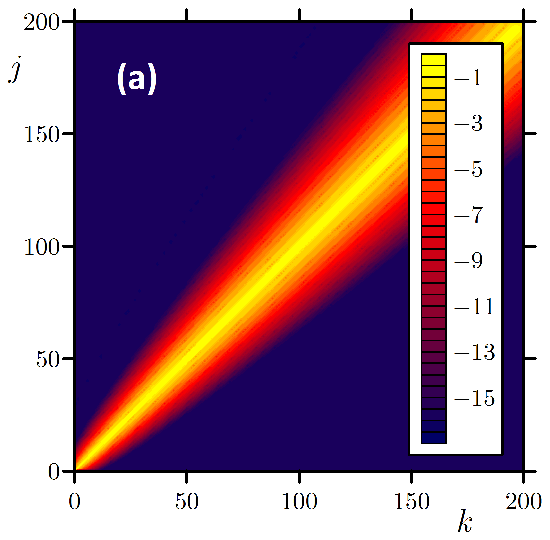}
\qquad\qquad
\includegraphics[width=0.35\textwidth]{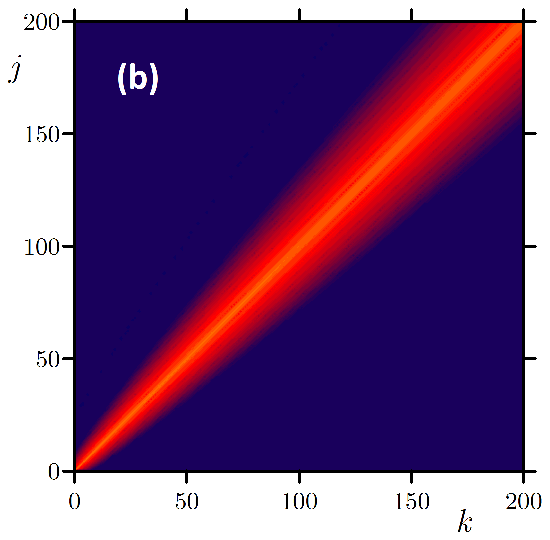}
}
\vspace{7pt}
\centerline{
\includegraphics[width=0.35\textwidth]{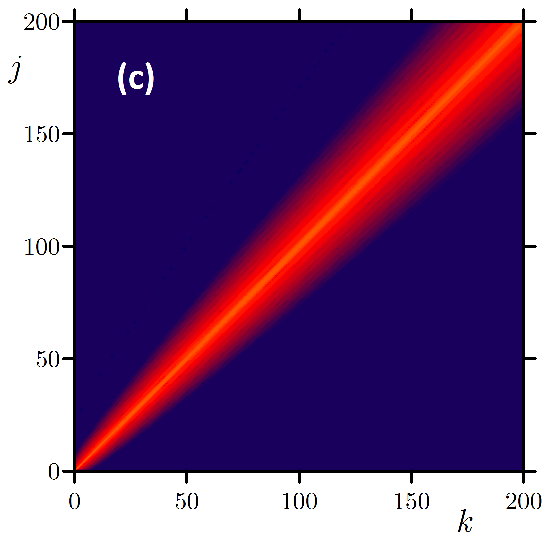}
\qquad\qquad
\includegraphics[width=0.35\textwidth]{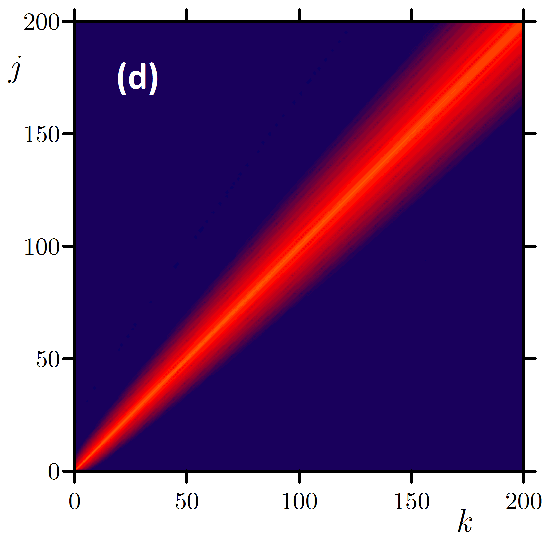}}
\caption{\label{fig9}
For equation system~(\ref{eq:109}) and (\ref{eq:110}), matrices
$\log_{10}|Q_{jk}|$ (panel \textit{a}),
$\log_{10}|(\mathbf{M}_1)_{jk}|$ (panel \textit{b}),
$\log_{10}|(\tau^{-1}\mathbf{M}_2)_{jk}|$ (panel \textit{c}),
$\log_{10}|(\tau^{-2}\mathbf{M}_3)_{jk}|$ (panel \textit{d})
are plotted versus indices $k$ and $j$ for $\tau=10^{-4}$, $\eta_0=-2$, $\gamma=1$, $J=10$, $\sigma^2=10$;
 $\tau_1=1.25\times10^{-6}$.
}
\end{figure*}

In the Fourier space, $w_\eta(\phi,t)=\frac{1}{2\pi}\big[1+\sum_{j=1}^\infty(a_j(t)e^{-ij\phi}+c.c.)\big]$, where ``$c.c.$'' stands for complex conjugate, and Eq.~(\ref{eq:104}) takes the following form:
\begin{align}
\dot{a}_j(\eta)&=j\left[2ia_j+\frac{i}{2}A_\eta(a_{j-1}+2a_j +a_{j+1})\right]
\nonumber\\
&
\quad
{}
-\sigma^2\left[\frac32j^2a_j+\bigg(j^2-\frac{j}{2}\bigg)a_{j-1}
 +\bigg(j^2+\frac{j}{2}\bigg)a_{j+1}
 +\frac{j(j-1)}{4}a_{j-2} +\frac{j(j+1)}{4}a_{j+2}\right]\;,
\label{eq:107}
\end{align}
where $A_\eta=\eta+Js+I(t)-1$; $a_0=1$ and $a_{-j}=a_j^\ast$, by definition. For a heterogeneous population with Lorentzian distribution $g(\eta)=\gamma/[\pi(\gamma^2+(\eta-\eta_0)^2)]$ with median $\eta_0$ and half-width at half maximum $\gamma$, one can employ the Residue theorem~\cite{Yakubovich-1969,Rabinovich-Trubetskov-1989,Crawford-1994} to derive from Eq.~(\ref{eq:107}) the infinite chain of dynamics equations for the Kuramoto--Daido order parameters~\cite{Daido-1996}
\[
z_j=\int\mathrm{d}\eta\,g(\eta)\,a_j(\eta) =\int\mathrm{d}\eta\,g(\eta)\int\mathrm{d}\phi\,w_\eta(\phi)e^{ij\phi}.
\]
One finds
\begin{align}
\dot{z}_j&=j\left[2iz_j+\frac{iA_\eta-\gamma}{2}(z_{j-1}+2z_j +z_{j+1})\right]
\nonumber\\
&
\qquad
{}
-\sigma^2\left[\frac32j^2z_j+\bigg(j^2-\frac{j}{2}\bigg)z_{j-1}
 +\bigg(j^2+\frac{j}{2}\bigg)z_{j+1}
 +\frac{j(j-1)}{4}z_{j-2} +\frac{j(j+1)}{4}z_{j+2}\right]
\label{eq:108}
\\
\qquad
&\equiv(\mathbf{L}\cdot\mathbf{z})_j+f_j(\mathbf{z},t)\,,
\label{eq:109}
\end{align}
where $\mathbf{z}=\{z_0,z_1,z_2,\dots\}$ (which does not contradict $z_0=1$),
\[
f_j=ij\frac{Jr+I(t)}{2}(z_{j-1}+2z_j +z_{j+1})\,,
\]
and everything else is collected in the linear part with constant coefficients.
Notice, the $Jr$-term of $\mathbf{f}$ is quadratic with respect to $\mathbf{z}$, as the firing rate $r(t)$ is a function of the system state. Namely, substituting the Fourier series of $w_\eta(\phi)$ into Eq.~(\ref{eq:105}), one finds the neuron firing rate~\cite{Montbrio-Pazo-Roxin-2015}
\begin{align}
r&=\frac{1}{\pi}(1-z_1-z_1^\ast+z_2+z_2^\ast-z_3-z_3^\ast+z_4+z_4^\ast+\dots)
\nonumber\\
&=\frac{\mathrm{Re}W}{\pi}\,,
\quad W=1-2z_1+2z_2-2z_3+2z_4+\dots\,.
\label{eq:110}
\end{align}

To summarize, the macroscopic dynamics of a recurrent synaptic network of QIFs with all-to-all coupling is governed by the system of equations~(\ref{eq:108}) and (\ref{eq:110}); its Eq.~(\ref{eq:ETD01}) form is given by (\ref{eq:109}) and (\ref{eq:110}).


\subsubsection{Direct numerical simulation challenges and ETD methods}
For perfect synchrony regimes, the distribution $w(\phi,t)$ is a travelling delta-function $\delta(\phi-\varphi(t))$, the Fourier spectrum of which is $z_j=e^{ij\varphi(t)}$ and does not decay with $j$. In physically meaningful set-ups the synchrony will be imperfect, but still high-synchrony regimes are of interest. The series $z_j$ for these regimes decay slowly and may require as many as $1000$--$2000$ elements for a reasonably accurate simulation of the dynamics of observables $z_1(t)$, $r(t)$, etc.~\cite{Goldobin-2021} In what follows, we estimate how much demanding for the numerical time stepsize are Eqs.~(\ref{eq:108}) for large $j$.


Typical explicit schemes which are stable for infinitesimally small time stepsize can become unstable for a larger stepsize $\tau_\ast$. For a perturbation $\delta z_j\sim C_ke^{ikj}$, where only $-\pi<k<\pi$ make sense as $j$ is integer, the reference value of $\tau_\ast$ can be estimated from the condition $|\delta\dot{z}_j\tau_\ast|/|\delta z_j|\sim 1$ (details vary from scheme to scheme, but the order of magnitude is $1$). For Eq.~(\ref{eq:108}) and large $j$, this condition reads:
\begin{align}
\left|j\left\{2i+(iA_\eta-\gamma)(1+\cos{k})\right\} -j\sigma^2\left\{j\left(\frac{3}{2}+2\cos{k}
+\frac{1}{2}\cos{2k}\right) +i\sin{k}+\frac{i}{2}\sin{2k}\right\}\right|\tau_\ast\sim 1\,.
\nonumber
\end{align}
Keeping only principal contributions and maximizing over $k$, one finds approximately
\[
\left(2j\sqrt{(1+A_\eta)^2+\gamma^2}+4j^2\sigma^2\right)\tau_\ast\sim 1\,.
\]
The Runge--Kutta--Merson and simplest predictor-corrector schemes yield stability threshold values of the same order of magnitude. For computationally challenging cases of $j\sim 1000$, the maximal admissible time stepsize for both conventional schemes becomes as small as $\tau_\ast\sim 10^{-7}/(\sigma^2+0.0005\sqrt{(1+A_\eta)^2+\gamma^2})$; however, this and even bigger number of modes is required for some high-synchrony regimes.

\begin{table*}[t]
  \centering
  \caption{
 CPU time for the direct numerical simulation of equation system~(\ref{eq:108}) and (\ref{eq:110}) is  averaged over 100 runs on the dimensionless time interval $t\in[0;5]$ for $\eta_0=-2$, $\gamma=1$, $J=10$, $\sigma^2=10$, $I(t)=0.3\sin{10t}$;
  ``preparation'': time for calculation of $\mathbf{Q}$ and $\mathbf{M}_n$ [sec], ``run'': time of simulation with a numerical scheme with pre-calculated matrixes $\mathbf{Q}$ and $\mathbf{M}_n$ [sec];
  stepsize $\tau$ was varied for the ETD schemes, and for the RKM- and PC-ones, the maximal values allowing for a numerically stable simulations were used;
  relative error $\varepsilon\equiv\sum_{j=1}^N|\delta z_j|/\sum_{j=1}^N|z_j|$.
  For the CPU specifications see Caption to Fig.~\ref{fig3}.
  }
\vspace{10pt}

%
\begin{tabular}{|p{2.1cm}|c|c|c|c|c|}
\hline
& \multirow{2}{*}{$N$} & \multirow{2}{*}{$\tau$} & \multicolumn{2}{c|}{simulation CPU time [sec]} & \multirow{2}{*}{$\varepsilon$} \\
\cline{4-5}
&  &  &\;preparation\;& run &  \\
\hline
\multirow{3}{2.1cm}{Run\-ge--Kut\-ta--Merson~method} &\; $200$ \;&\,$2.3\times 10^{-6}$\,& \multicolumn{2}{c|}{$16.1$} &  \\
 & $300$ & $10^{-6}$ & \multicolumn{2}{c|}{$53.0$} &  \\
 & $400$ & $0.5\times 10^{-6}$ & \multicolumn{2}{c|}{$139$} &  \\
\hline
\multirow{3}{2.1cm}{predictor-corrector~(\ref{eq:CH05ad})} & $200$ & $1.2\times 10^{-6}$ & \multicolumn{2}{c|}{$13.8$} &\,$6.1\times10^{-11}$\,\\
 & $300$ & $5.5\times 10^{-7}$ & \multicolumn{2}{c|}{$49.9$} & $1.4\times10^{-11}$ \\
 & $400$ & $3.1\times 10^{-7}$ & \multicolumn{2}{c|}{$106$} & $4.9\times10^{-10}$ \\
\hline
\multirow{6}{*}{ETD2RK} &  \multirow{2}{*}{200} & $0.0005$ & $0.752$ & $0.936$ & $5.1\times10^{-8}$ \\
& & $0.005$ & $7.51$ & $0.149$ & $5.0\times10^{-6}$ \\
& \multirow{2}{*}{300} & $0.0005$ & $4.23$ & $1.96$ & $5.1\times10^{-8}$ \\
& & $0.005$ & $42.1$ & $0.300$ & $5.0\times10^{-6}$ \\
& \multirow{2}{*}{400} & $0.0005$ & $12.2$ & $3.58$ & $5.1\times10^{-8}$ \\
& & $0.005$ & $119$ & $0.656$ & $5.0\times10^{-6}$ \\
\hline
\multirow{6}{*}{ETD3RK} &  \multirow{2}{*}{200} & $0.0005$ & $1.25$ & $2.12$ & $1.1\times10^{-10}$ \\
& & $0.005$ & $12.5$ & $0.345$ & $1.0\times10^{-7}$ \\
& \multirow{2}{*}{300} & $0.0005$ & $7.05$ & $5.56$ & $1.9\times10^{-10}$ \\
& & $0.005$ & $70.1$ & $1.21$ & $1.0\times10^{-7}$ \\
& \multirow{2}{*}{400} & $0.0005$ & $20.3$ & $17.6$ & $4.3\times10^{-10}$ \\
& & $0.005$ & $199$ & $3.23$ & $1.0\times10^{-7}$ \\
\hline
\multirow{6}{*}{ETD4RK} &  \multirow{2}{*}{200} & $0.0005$ & $1.25$ & $2.59$ & $5.6\times10^{-11}$ \\
& & $0.005$ & $12.5$ & $0.397$ & $4.2\times10^{-10}$ \\
& \multirow{2}{*}{300} & $0.0005$ & $7.05$ & $6.02$ & $1.5\times10^{-10}$ \\
& & $0.005$ & $70.1$ & $1.11$ & $4.5\times10^{-10}$ \\
& \multirow{2}{*}{400} & $0.0005$ & $20.3$ & $17.5$ & $4.5\times10^{-10}$ \\
& & $0.005$ & $199$ & $3.36$ & $7.9\times10^{-10}$ \\
\hline
\end{tabular}
  \label{tab1}
\end{table*}

For the application of an ETD method for this problem, one practically cannot make the matrix $\mathbf{L}$ diagonal or calculate matrix $\exp(\mathbf{L}\tau)$ analytically in some other way. On the other hand, one can calculate matrixes $\mathbf{Q}$, $\mathbf{M}_n$, etc.\ numerically with the Runge--Kutta--Merson or simplest predictor-corrector schemes. In Fig.~\ref{fig8}, one can see the results of the direct numerical simulation of system~(\ref{eq:108}) and (\ref{eq:110}) with ETD methods for $N=200$ modes $z_j$. The reference ``accurate'' solution is calculated with the predictor-corrector scheme~(\ref{eq:CH05ad}) and small time step $\tau=1.25\times 10^{-7}$.~\footnote{The stepsize around $\tau\sim10^{-7}$ was found to provide the best accuracy with double precision calculations for the given set of parameter values. Here, the accumulation of the relative machine error $\sim10^{-15}$ at each time step sums up over the unit time to a value much smaller than $\sim10^{-15}/10^{-7}=10^{-8}$ because of the fast decay of the leading perturbation modes; as a result, the net numerical error rate per the unit time is well below $10^{-12}$.}
The ETD scheme stepsize is set to a large value $\tau=0.025$ in order to illustrate the accuracy order of three employed ETD schemes. In Fig.~\ref{fig9}, the structure of numerically calculated matrices $\mathbf{Q}$, $\mathbf{M}_n$ is presented for $\tau=10^{-4}$.

In Table~\ref{tab1}, the performance of three EDT schemes is compared against the background of the conventional Runge--Kutta--Merson and predictor-corrector~(\ref{eq:CH05ad}) schemes. For typical studies, as many modes as $N=200$--$400$ is sufficient for accurate resolution of fine elements of phase diagrams. For short-term simulations, optimal stepsize was found to be around $\tau=0.0005$ and ETD schemes provide the performance gain in range $3$--$10$; for long-term simulations, where one can neglect the time for calculation of the scheme coefficients, a bigger stepsize $\tau=0.005$ allows one to have decent accuracy alongside with the performance gain up to two orders of magnitude.

\section{Conclusion \label{sec:Concl}}

We have suggested a practical approach for the implementation of the exponential time differencing schemes of the Runge--Kutta type~\cite{Cox-Matthews-2002} for numerical simulation of stiff systems with nondiagonal linear part. In this approach, the numerical schemes are written in terms of the solutions of auxiliary problems which are preliminary integrated numerically by means of an explicit basic method with a very small time step but over a short interval of time---one step of the ETD scheme. As a result, one can not only use the ETD methods for the systems where an analytical calculation of the exponential of a nondiagonal linear operator is infeasible, but also employ the same program code for different equation systems, the only part of code which is subject to change is the subroutine with the equations to be simulated.

The employment of this approach for Cahn--Hilliard equation~(\ref{eq:CH01}) and sixth-order spatial derivative Matthews--Cox equation~(\ref{eq:MC01}) demonstrated that one can have a performance gain by two orders of magnitude for the former and by three orders for the latter. Moreover, for partial differential equations, the performance gain grows as one increases the required accuracy of spatial resolution; the basic gain for an equation with the highest-order spatial derivative $\partial^m/\partial x^m$ is $\propto h_x^{m-2}$. For nonlarge $m$, one can introduce code optimization accounting for the sparseness of the matrices $\mathbf{Q}$ and $\mathbf{M}_n$, which gives some further performance acceleration on top of the basic one (this acceleration is by far not that large as the basic gain). For the simulation of a long-term evolution, as required for complex spatiotemporal dynamics or Anderson localization phenomena, the maximal performance gain is achieved with the 4th-order Runge--Kutta type ETD scheme; for a short-time simulation, the usage of the 2nd-order ETD scheme is somewhat more efficient.

The employment of the approach for spectral methods is illustrated with 1D Fokker--Planck equation governing the macroscopic dynamics of populations of quadratic integrate-and-fire neurons with endogenic noise. Here, one can have a controllable accuracy and achieve the performance gain up to one order of magnitude for time-independent regimes, where short-term simulations are sufficient, and several orders of magnitude for time-dependent regimes where long-term simulations are demanded. Our approach is also efficient for numerical simulation of other classical systems such as populations of coupled active rotators~\cite{Sakaguchi-Shinomoto-Kuramoto-1988,Klinshov-etal-2021}.

The suggested approach can be extended into multiple dimensions in a straightforward way~\cite{Permyakova-Goldobin-2020}. However, the basic version of the code for a partial differential equation, say, in 2d requires matrices $\mathbf{Q}$ and $\mathbf{M}_n$ of size $[N\times N]$ with $N=N_x\times N_y$, where $N_x$ and $N_y$ are the number of nodes (or modes) in the $x$- and $y$-directions. With the optimization for the sparseness of $\mathbf{Q}$ and $\mathbf{M}_n$, one can store a reduced amount $[N\times N_\mathrm{spr}]$ of elements of each matrix, where $N_\mathrm{spr}=2m\tau\min\left\{N_y|c_x|\ln[h_x/(err\,\tau)],N_x|c_y|\ln[h_y/(err\,\tau)]\right\}$ for the highest-order spatial derivative term $c_x\partial^mu/\partial x^m+c_y\partial^mu/\partial y^m$ and the acceptable relative error per a unit time $err$. For high-precision simulations in multiple dimensions the usage of the ETD schemes becomes memory-demanding. Elevated memory requirements and the fact, that the bulk of the computation process is the multiplication of large matrices and vectors, make these high-performance simulation methods naturally suitable for parallel- and super-computing.

\begin{figure*}[t]
\centerline{
\includegraphics[width=0.42\textwidth]{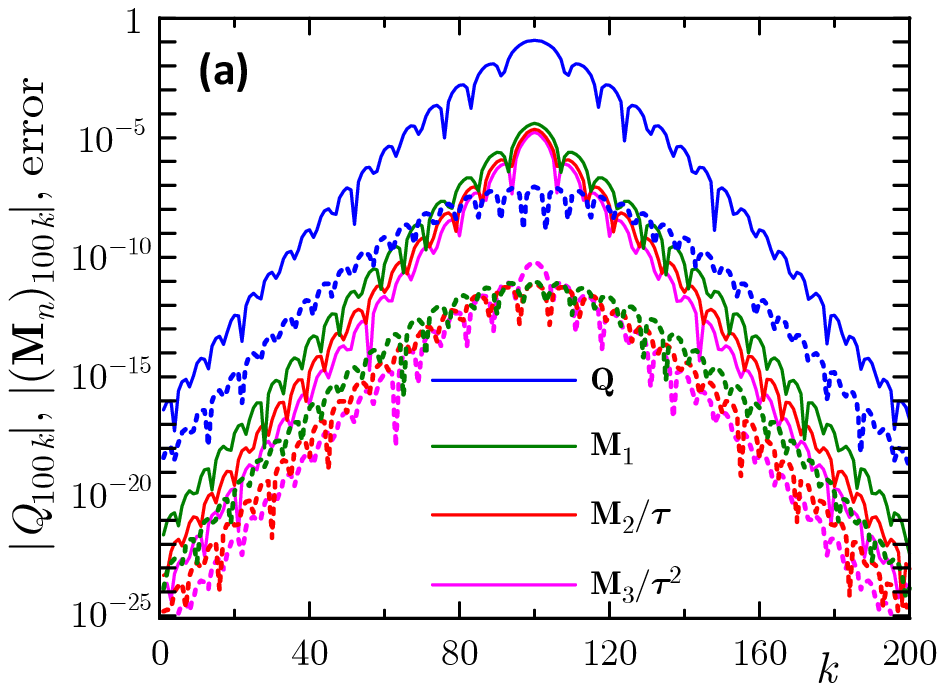}
\qquad\qquad
\includegraphics[width=0.42\textwidth]{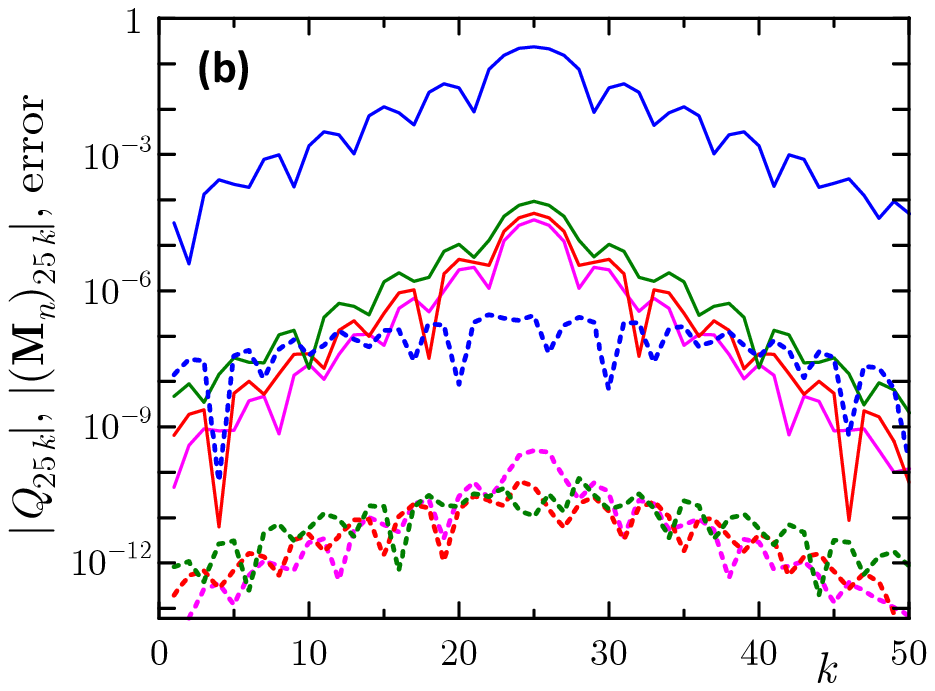}}
\caption{\label{fig10}
The strong accuracy of calculation of $Q_{jk}$, $(\mathbf{M}_1)_{jk}$, $(\tau^{-1}\mathbf{M}_2)_{jk}$, $(\tau^{-2}\mathbf{M}_3)_{jk}$ is demonstrated with the values of coefficients for CHE~(\ref{eq:CH01}) and $j=100$ (panel~\textit{a}: $\tau=2.5\times10^{-4}$, $N=200$, i.e.\ $h_x=0.05$, other parameter values are as in Figs.~\ref{fig2} and \ref{fig3}) and for MCE~(\ref{eq:MC01}) and $j=25$ (panel~\textit{b}: $\tau=3.2\times10^{-4}$, $N=50$, i.e.\ $h_x=0.2$, other parameter values are as in Fig.~\ref{fig6}). Solid lines: ``exact'' values; dotted lines: error; the color coding is identical in panels \textit{a} and \textit{b}. The ``exact'' coefficients are calculated with Taylor series of length $M=2088$ and $384$ decimal digits of intermediate computations (\textit{a}) and $M=1044$ and $192$ decimal digits (\textit{b}).
}
\end{figure*}
\begin{figure*}[t]
\centerline{
\includegraphics[width=0.52\textwidth]{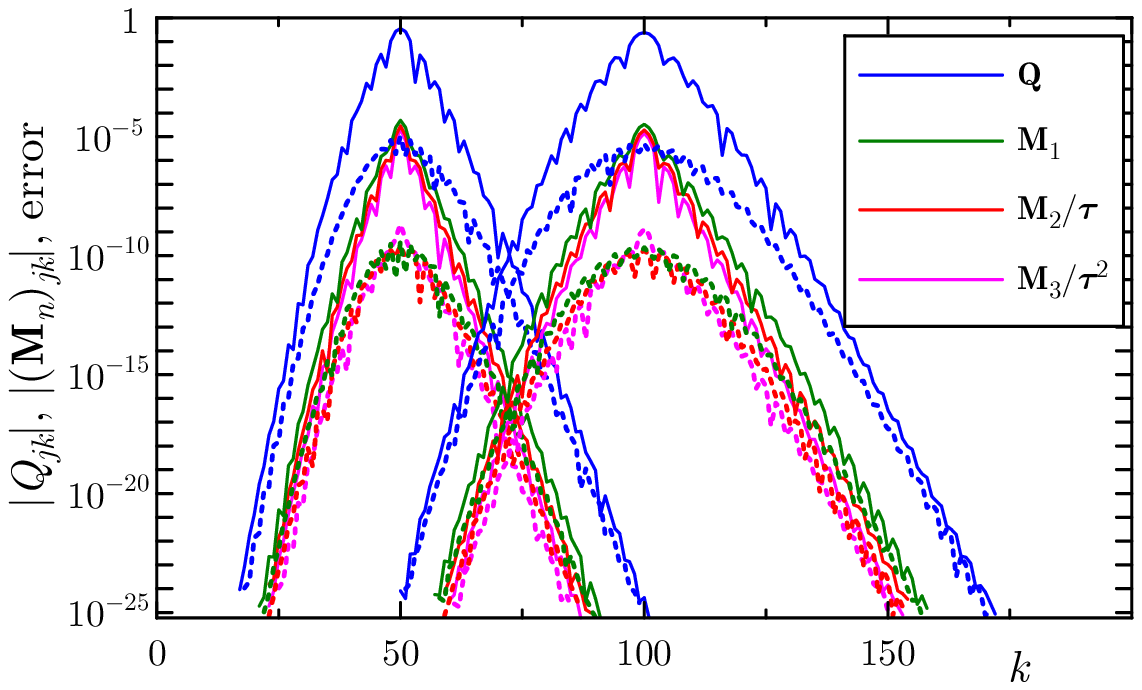}
}
\caption{\label{fig11}
For Fokker--Planck equation~(\ref{eq:108}), ``exact'' coefficients
$Q_{jk}$, $(\mathbf{M}_1)_{jk}$, $(\tau^{-1}\mathbf{M}_2)_{jk}$, $(\tau^{-2}\mathbf{M}_3)_{jk}$
are plotted versus index $k$ with the solid lines for $j=50$ (the left set of curves) and $100$ (the right set of curves); the error of the coefficients calculated with the proposed procedure is plotted with the dotted curves.
The ``exact'' coefficients are calculated with Taylor series of length $M=522$ and $96$ decimal digits.
See Caption to Fig.~\ref{fig9} for the parameter values.
}
\end{figure*}

\appendix
\section{Strong accuracy of computation of $\mathbf{Q}$ and $\mathbf{M}_n$}\label{sec:app}
In this section we evaluate the ``strong'' accuracy of the approximate calculation of matrices $\mathbf{Q}$ and $\mathbf{M}_n$ with the suggested algorithm. Below we will also discuss why one distinguishes the ``strong'' accuracy of calculation of these coefficients and the accuracy of the numerical solution $\mathbf{u}(t)$, which is the ``soft'' accuracy of the numerical simulation of dynamical system~(\ref{eq:ETD01}) and often higher; significantly higher in the examples we considered in this paper.

Matrices $\mathbf{Q}=e^{\tau\mathbf{L}}$ and $\mathbf{M}_n$ given by Eqs.~(\ref{eq:ETDM13}) are determined by their series:
\begin{align}
\mathbf{Q}&=e^{\tau\mathbf{L}}=\mathbf{I}+\tau\mathbf{L}+\frac{\tau^2}{2!}\mathbf{L}^2+\frac{\tau^3}{3!}\mathbf{L}^3+\dots\,,
\label{eq:app101}\\
\mathbf{M}_n&=\frac{\tau^n}{n}\mathbf{I}+\frac{\tau^{n+1}\mathbf{L}}{n(n+1)}+\frac{\tau^{n+2}\mathbf{L}^2}{n(n+1)(n+2)}+\dots +\frac{\tau^{n+m}\mathbf{L}^{n+m}}{n(n+1)\dots(n+m)}+\dots\,.
\label{eq:app102}
\end{align}
Notice, the set of formulas~(\ref{eq:ETDM13}) is just a short way to write series~(\ref{eq:app102}) and they are well defined for the case where matrix $\mathbf{L}$ possesses zero eigenvalues (which are inherent to conservation laws in physical systems and some bifurcation points) and $\mathbf{L}^{-1}$ diverges.

A mathematically idealistic way to compute these matrices is related to their diagonalization or decomposition into the basis of eigenvectors $\mathbf{y}_j$ ($j=1,2,...,N$) of $\mathbf{L}$, defined by $\mathbf{L}\cdot\mathbf{y}_j=\lambda_j\mathbf{y}_j$ with the corresponding Hermitian adjoint problems $\mathbf{v}_j^T\cdot\mathbf{L}=\lambda_j\mathbf{v}_j^T$, where the superscript ``$T$'' indicates a transposed vector/matrix and the orthonormalization condition is fulfilled: $\mathbf{v}_j^T\cdot\mathbf{y}_m=\delta_{jm}$. With a solved eigenvalue problem, one can substitute $\mathbf{L}=\sum_{j=1}^{N}\lambda_j\mathbf{y}_j\mathbf{v}_j^T$ into series (\ref{eq:app101}) and (\ref{eq:app102}) and calculate
\begin{align}
\mathbf{Q}&=\sum_{j=1}^{N}e^{\lambda_j\tau}\mathbf{y}_j\mathbf{v}_j^T\,,
\label{eq:app103}\\
\mathbf{M}_n&=\sum_{j=1}^{N}\frac{(n-1)!}{(\lambda_j)^n}\left(e^{\lambda_j\tau} -\sum_{m=0}^{n-1}\frac{(\lambda_j\tau)^m}{m!}\right)\mathbf{y}_j\mathbf{v}_j^T\,.
\label{eq:app104}
\end{align}
The $j$-th summand in $\mathbf{M}_n$ is again correctly defined also for $\lambda_j=0$: one takes the limit $\lambda_j\to0$ and finds
$(\tau^n/n)\mathbf{y}_j\mathbf{v}_j^T$.
In practice, the exact solution of the eigenvalue problem is not always feasible and approaches for high-precision approximate calculations are available in the literature (e.g., employing the projection into Krylov subspaces~\cite{Hochbruck-Lubich-Selhofer-1998,Tokman-2006}).

Our task in this section is to calculate the ``exact'' matrices $\mathbf{Q}$ and $\mathbf{M}_n$ with a controlled accuracy (arbitrary small error) and examine the results of their calculations with the algorithm of this paper. Therefore, we perform a calculation with Taylor series~(\ref{eq:app101}) and (\ref{eq:app102}) with as many elements and decimal digits as needed to have the results with the standard double precision accuracy---not less than $16$ decimal digits.
The Taylor series of the exponentials of $(\lambda_j\tau)$ in Eqs.~(\ref{eq:app103}) and (\ref{eq:app104}) with respect to $\tau$ are equivalent to series (\ref{eq:app101}) and (\ref{eq:app102}). Since the exponential does converge for all argument values, the series we have to compute always converges even thought a gigantic number of elements can be required.

The most operation-number efficient way of computing $\mathbf{Q}$ employs the formula (\ref{eq:app101}) in a form $\mathbf{Q}=\mathbf{I}+\tau\mathbf{L}\cdot\Big(\mathbf{I} +\frac{\tau}{2}\mathbf{L}\cdot\Big(\mathbf{I} +\frac{\tau}{3}\mathbf{L}\cdot\Big(\mathbf{I}+\dots\Big)\Big)\Big)$. Accordingly to this form, for a Taylor series truncated after the $M$-th element, one can write $\mathbf{Q}_{M}=\mathbf{I}$ and iteratively compute
\begin{equation}
\mathbf{Q}_{m-1}=\mathbf{I} +\frac{\tau}{m}\mathbf{L}\cdot\mathbf{Q}_{m}\,,
\label{eq:app105}
\end{equation}
descending from $m=M$ to $m=1$. This iterative procedure yields all required matrices,
\begin{equation}
\mathbf{Q}=\mathbf{Q}_0\,,
\qquad
\mathbf{M}_n=\frac{\tau^n}{n}\mathbf{Q}_n\,,
\label{eq:app106}
\end{equation}
and provides the slowest accumulation rate for the machine cancellation error.

The core purpose of the ETD methods is to iterate a numerical scheme with time stepsize $\tau$ for which $\lambda_\mathrm{fast}\tau\gg1$, where $\lambda_\mathrm{fast}$ is the eigenvalue of the fastest decaying or oscillating mode. Hence, the practically relevant cases require the calculation of Taylor series with as many elements as needed to resolve $e^{\lambda_\mathrm{fast}\tau}$ [see Eqs.~(\ref{eq:app103}) and (\ref{eq:app104}), which are equivalent to (\ref{eq:app101}) and (\ref{eq:app102})] for a large $\lambda_\mathrm{fast}\tau$ or at least avoid its numerical explosion. To avoid the numerical instability of the fastest mode caused by the Taylor series truncation, one needs the order of magnitude of the truncated terms for this mode $\sim|\lambda_\mathrm{fast}\tau|^M/M!<1$. With the Stirling's approximation $M!\approx\sqrt{2\pi M}(M/e)^M$, one can estimate the minimal required truncation length $M_\ast$:
\[
\frac{e|\lambda_\mathrm{fast}|\tau}{M_\ast}\frac{1}{(2\pi M_\ast)^\frac{1}{2M_\ast}}\sim 1\,.
\]
Factor $(2\pi M_\ast)^\frac{1}{2M_\ast}$ tends to 1 for large $M_\ast$ and is always bigger than $1$;
therefore, we can safely adopt
\begin{equation}
M_\ast=e|\lambda_\mathrm{fast}|\tau\,.
\label{eq:app107}
\end{equation}

For CHE~(\ref{eq:CH01}), the short wavelength modes
$\delta u(x,t)\propto c_\kappa(t)e^{i\kappa x}$ with the wavenumber $\kappa\gg1$ and amplitude $c_\kappa(t)$ are the fastest decaying ones. Their decay rate is dominated by the last term in~(\ref{eq:CH03}) and, substituting the perturbation $c_\kappa(t)e^{i\kappa x}$ into~(\ref{eq:CH03}), one finds: $\lambda_\kappa c_\kappa(t)e^{i\kappa x}=\big[\dots-\big(e^{i\kappa h_x/2}-e^{-i\kappa h_x/2}\big)^4/h_x^4\big]c_\kappa(t)e^{i\kappa x}=\big[\dots-(16/h_x^4)\sin^4(\kappa h_x/2)\big]c_\kappa(t)e^{i\kappa x}$. Hence, $\lambda_\mathrm{fast}\approx-16/h_x^4$\, and Eq.~(\ref{eq:app107}) yields
\[
M_\ast=16e\tau/h_x^4\,.
\]
For the example in Fig.~\ref{fig2} with $\tau=2.5\times10^{-4}$ and the same values of other parameters, we find $M_\ast=1740$ and take truncation length $M=1.2M_\ast=2088$ to be on a safe side with the series convergence. We use the Maple analytical calculation package for computations with a variable number of decimal digits to handle the problem of cancelation error accumulation. For $256$ decimal digits the first $20$ digits of $Q_{jk}$ are still affected by the variation of the mantissa length (cancellation error). We conduct the calculations with $384$ digits and find the first 20 digits of the results to be the same as with $512$ for $M=1.2M_\ast$ and elongated series $M=1.4M_\ast=2436$. Therefore, we trust the first $20$ digits of the ``exact'' calculated matrices $\mathbf{Q}$ and $\mathbf{M}_n$. In Fig.~\ref{fig10}\textit{a}, one can see that the relative error of the significant parts of $Q_{jk}$ and $(\mathbf{M}_n)_{jk}$ is somewhat smaller than $10^{-6}$.

With the approach of this paper the fine dynamics of the fastest decaying modes is resolved accurately enough to avoid the numerical instability but not with a high precision. The corresponding contributions into $\mathbf{Q}$ and $\mathbf{M}_n$ are not very accurate. These contributions result in the ``strong'' accuracy of calculation of the matrices $\lesssim10^{-6}$ for the example in Figs.~\ref{fig10}\textit{a} and \ref{fig2}. However, the inaccuracy of the fastest decaying perturbations makes an exponentially small contribution into the numerically simulated system state; the accuracy of numerical solutions is practically important and is much higher as one can see in Figs.~\ref{fig3} and \ref{fig4}. The latter can be referred to as the ``soft'' accuracy of the numerical scheme.

Similarly, for MCE~(\ref{eq:MC01}), the short wavelength modes
$\delta u(x,t)\propto c_\kappa(t)e^{i\kappa x}$ are the fastest decaying ones. Their decay rate is dominated by the last term in~(\ref{eq:MC02a}) and one finds: $\lambda_\kappa c_\kappa(t)e^{i\kappa x}=\big[\dots+\big(e^{i\kappa h_x/2}-e^{-i\kappa h_x/2}\big)^6/h_x^6\big]c_\kappa(t)e^{i\kappa x}=\big[\dots-(64/h_x^6)\sin^6(\kappa h_x/2)\big]c_\kappa(t)e^{i\kappa x}$. Hence, $\lambda_\mathrm{fast}\approx-64/h_x^6$\, and Eq.~(\ref{eq:app107}) yields
\[
M_\ast=64e\tau/h_x^6\,.
\]
For the example in Fig.~\ref{fig6} with $\tau=3.2\times10^{-4}$ and the same values of other parameters, we find $M_\ast=870$ and take truncation length $M=1.2M_\ast=1044$. For calculations with $128$ decimal digits the first $20$ digits of $Q_{jk}$ are still affected by the cancellation error. We conduct the calculations with $192$ digits and find the first 20 digits of the results to be the same as with $256$ for $M=1.2M_\ast$ and $M=1.4M_\ast=1218$. Therefore, we trust the first $20$ digits of the ``exact'' calculated matrices $\mathbf{Q}$ and $\mathbf{M}_n$. In Fig.~\ref{fig10}\textit{b}, one can see that the relative error of the significant parts of $Q_{jk}$ and $(\mathbf{M}_n)_{jk}$ is $\lesssim10^{-6}$. The ``soft'' accuracy of the numerical scheme is higher, as one can see in Figs.~\ref{fig6} and \ref{fig7}.

For Fokker--Planck equation~(\ref{eq:108}), the fastest decaying modes
$\delta z_j(t)\propto c_\kappa(t)e^{i\kappa j}$ are localized at largest $j\approx N$.  Their decay rate is dominated by the $\sigma^2$-term in~(\ref{eq:108}) and one finds: $\lambda_\kappa c_\kappa(t)e^{i\kappa j}=\big[\dots-\sigma^2j^2\big(e^{i\kappa/2}+e^{-i\kappa/2}\big)^4/4\big]c_\kappa(t)e^{i\kappa j}=\big[\dots-4j^2\sigma^2\cos^4(\kappa/2)\big]c_\kappa(t)e^{i\kappa j}$. Hence, $\lambda_\mathrm{fast}\approx-4N^2\sigma^2$\, and Eq.~(\ref{eq:app107}) yields
\[
M_\ast=4e\sigma^2N^2\tau\,.
\]
For the example in Fig.~\ref{fig9} with $\tau=10^{-4}$ and the same values of other parameters, we find $M_\ast=435$ and take truncation length $M=1.2M_\ast=522$. For calculations with $64$ decimal digits the first $20$ digits of $Q_{jk}$ are still affected by the cancellation error. We conduct the calculations with $96$ digits and find the first 20 digits of the results to be the same as with $128$ for $M=1.2M_\ast$ and $M=1.4M_\ast=609$. Therefore, we trust the first $20$ digits of the ``exact'' calculated matrices $\mathbf{Q}$ and $\mathbf{M}_n$. In Fig.~\ref{fig11}, one can see the relative error of the significant parts of $Q_{jk}$ and $(\mathbf{M}_n)_{jk}$ is somewhat smaller than $10^{-5}$. The ``soft'' accuracy of the numerical scheme is still better, as one can see in Fig.~\ref{fig8} and Table~\ref{tab1}.

\section*{CRediT authorship contribution statement}

{\bf Evelina V.\ Permyakova:}
Conceptualization (supporting), Methodology (equal), Software (equal), Validation (equal), Formal analysis (equal), Visualization (supporting), Writing -- original draft (supporting), Writing -- review \& editing (supporting).
{\bf Denis S.\ Goldobin:}
Conceptualization (lead), Methodology (equal), Software (equal), Validation (equal), Formal analysis (equal), Visualization (lead), Writing -- original draft (lead), Writing -- review \& editing (lead).

\section*{Declaration of competing interest}
The authors declare no competing interests to disclose.

\section*{Data availability statement}

The data that supports the findings of this study are available within the article in the graphic form. The data sheets for the graphs and the program codes in FORTRAN are available on request from the authors.

\section*{Acknowledgements}
The work was carried out as part of a major scientific project (Agreement no.\ 075-15-2024-535 by April 23, 2024).

\bibliographystyle{elsarticle-num}
\providecommand{\noopsort}[1]{}\providecommand{\singleletter}[1]{#1}%

\end{document}